\documentclass[11pt]{amsart}
\usepackage{anyfontsize}
\usepackage{tikz-cd}
\usepackage[headings]{fullpage}
\usepackage{amsmath}
\usepackage{amsfonts}
\usepackage{latexsym}
\usepackage{amssymb}
\usepackage{placeins}

\newtheorem{proposition}{Proposition}[section]
\newtheorem{theorem}{Theorem}[section]

\newtheorem{lemma}{Lemma}[section]
\newtheorem{definition}{Definition}[section]
\newtheorem{corollary}{Corollary}[section]

\def\char{\sf char}

\def\min{\mathop{\sf min}}

\def\homo{\mathop{\sf Hom}}
\def\ker{\mathop{\sf Ker}}

\def\ann{\mathop{\sf ann}}

\def\ass{\mathop{\sf Ass}}
\def\supp{\mathop{\sf supp}}

\def\t{{\sf T}}

\def\min{{\sf min}}

\def\sig{\sigma, \sigma^{-1}}

\def\B ellk{\mathcal{B}_{\ell,k}}

\def\C{\mathbb{C}}

\def\N{\mathbb{N}}
\def\R{\mathbb{R}}
\def\Z{\mathbb{Z}}

\def\gal{\mathop{\sf Gal}}
\def\aut{\mathop{\sf Aut}}
\def\gl{\mathop{\sf GL_n(\mathbb{Z})}}
\def\gltwo{\mathop{\sf GL_2(\mathbb{Z})}}
\def\pgl{\mathop{\sf PGL}}

\newcommand{\diag}{\mathsf{diag}}

\newcommand{\diff}{[\sigma_1, \sigma_1^{-1}, \ldots ,\sigma_n, \sigma_n^{-1}]}

\newcommand{\F}{\mathcal{F}}

\begin{document}

\title[Coarsest Lattice Determining a System]{The Coarsest Lattice That Determines a Discrete Multidimensional System}

\author{Debasattam Pal} 
\address{Department of Electrical Engineering, IIT Bombay} \email{debasattam@ee.iitb.ac.in}
\author{Shiva Shankar} 
\address{Department of Electrical Engineering, IIT Bombay} \email{shunyashankar@gmail.com } 
\maketitle

\begin{abstract} A discrete multidimensional system is the set of solutions to a system of linear partial difference equations defined on the lattice $\Z^n$. This paper shows that it is determined by a unique coarsest sublattice, in the sense that the solutions of the system on this sublattice determine the solutions on $\Z^n$; it is therefore the correct domain of definition of the discrete system. 
In turn,  the defining sublattice is determined by a Galois group of symmetries that  leave invariant the equations defining the  system. These results find application in understanding properties of the system such as controllability and autonomy, and in its order reduction. 
\end{abstract}

\vspace{3mm}


\section{introduction}
In this paper we study the process of {\it contracting} the trajectories of a discrete multidimensional system - in short, an $n$-D system - on the lattice $\Z^n$ to a sublattice $\mathbb{S}$, the reverse process of {\it extending} the trajectories of a system on $\mathbb{S}$ to the entire lattice, as well as the composite process of contraction followed by extension. Of special interest are systems which are invariant under this cycle of contraction and extension. Contracting such an invariant  system to $\mathbb{S}$ is analogous to the process of restricting the solutions of a partial differential equation to an invariant subset, or the flow of a vector field to an invariant manifold. Furthermore, being an extension, its trajectories on $\Z^n$ are constructed from those on $\mathbb{S}$ without the imposition of any further laws, or restrictions. The study of the properties of the extended system is thus reduced to its study on the sublattice $\mathbb{S}$; its domain of definition should be considered to be $\mathbb{S}$ rather than $\mathbb{Z}^n$. We show that an $n$-D system is defined by its trajectories on a unique coarsest sublattice, which then should be considered its correct domain of definition. 
Furthermore, we show that this  sublattice is characterised by a group of symmetries that leave invariant the equations defining the system.  

We now expand on the above points. An $n$-D system is a collection of trajectories, i.e. functions on $\Z^n$, each one of which obeys the laws of the system (the next section contains precise definitions). A law is a difference equation, an element in the Laurent polynomial ring $A = \C\diff$ of difference operators on $\Z^n$, where $\sigma_i$ denotes shift in the $i$-th direction. Such a  law relates the value of a function at a point of $\Z^n$ to its values at other points of $\Z^n$, and only those functions that satisfy the law qualify to be trajectories.  For example, the law $1 + \sigma - \sigma^2$ describes a 1-D system on $\Z$, each of whose trajectories  $f$ satisfies the equation $f(x) + f(x+1) - f(x+2) = 0$. 
When the trajectory is defined by some $k$ attributes of the system, then a law is an element of $A^k$. For instance, the law $(1 + \sigma - \sigma^2, ~\sigma)$ on $\Z$ relates the attributes $(f_1, f_2)$ of the system that it defines by the equation $f_1(x) + f_1(x+1) - f_1(x+2) + f_2(x+1) = 0$. For this reason, we use the term solution interchangeably with trajectory. 

In these terms, the phenomenon we investigate in this paper can be stated as follows: it might be that the value of a trajectory at some point is related to its values at points only in some subset of the lattice. Specifically, let $\mathbb{S}$ be a sublattice of $\Z^n$ as above, and suppose that the value of a trajectory at a point is related only to its values at other points in the coset of $\mathbb{S}$ in $\Z^n$ that it belongs to.
For instance, consider the system defined by $1+ \sigma^2 - \sigma^4 $ on $\Z$. Then, a trajectory $f$ satisfies $f(x) + f(x+2) - f(x+4) = 0$, and the value of $f$ at a point of $\Z$ depends only on its values at points in the coset of $2\Z$ in $\Z$ that it lies in. 
What can we say about such a system? How is its behaviour on $\Z^n$ related to its behaviour on $\mathbb{S}$? How do we detect the presence of such sublattices? In this paper we provide answers through an analysis of the symmetries of the system. Thus, let $A_S$ be the ring of difference operators on $\mathbb{S}$. If $\mathbb{S}$ is a sublattice of full rank, then the group $G$ of automorphisms of $A$ which leave $A_S$ fixed is a finite group, and a system can be reconstructed from its contraction to $\mathbb{S}$ if and only if the equations defining the system are left invariant by $G$. We also study the case of degenerate sublattices by reducing it to the full rank case. These answers provide a complete solution to our problem, which also highlights the fundamental role of symmetries in the study of $n$-D systems.

The complexity of a dataset is usually defined to be the minimum possible order of a system that generates it. For 1-D systems, Willems formalizes it by defining a notion of `memory' of the system \cite{jcw}, but the multidimensional version is more involved \cite{debp, deb, mp}. The results of our paper suggest an alternative, namely the coarser the sublattice from which a system can be reconstructed, the lower is its complexity. In the example above, the system $1 + \sigma^2 - \sigma^4$ is defined by a fourth order operator, but on $2\Z$, it is defined by $1 + \tau - \tau^2$, where $\tau$ is the shift operator on $2\Z$. As the system can be reconstructed from $2\Z$ without the imposition of any further laws, the true measure of its complexity should be determined by its behaviour on the sublattice. Thus, our article presents a novel notion of order reduction for $n$-D systems. This has immediate implications for computational solutions to partial difference equations using computer algebra packages that implement Gr\"{o}bner bases algorithms as in \cite{zo}. A detailed study of these implications will be pursued elsewhere.

We start our paper by studying separately each of the two processes of contraction and extension, independent of any requirement of invariance. Contracting the trajectories of any $n$-D system to a sublattice yields a system on the sublattice, and there are usually many systems on $\Z^n$ that contract to a given system on the sublattice. The invariant system of the above paragraphs is one of this set, its unique maximal element, but the others are important as well. The group $G$ permutes these systems, and the invariant one is the unique fixed point of this action. Similarly, the process of extension is important in its own right, and is the algebraic analogue of a fractional partial differential equation. For example, the extension from the sublattice $2\Z$ to $\Z$ is equivalent to adjoining the square roots of the shift operator.



The paper is organised as follows. We first consider the case of diagonal sublattices of $\Z^n$, namely sublattices of the type $\{(x_1, \ldots, x_n)~|~ x_i = 0, \pm d_i, \pm 2d_i, \ldots \}$. This is because all our calculations are then straightforward and transparent. Besides, a lot of the literature in the subject, including definitions such as `degree of autonomy', confine themselves to diagonal sublattices \cite{wro2,dp}. We set up notation in Section 2; in Section 3 we introduce the notions of contraction and extension. In Section 4, we use the isomorphism provided by the Smith Canonical Form to move a general sublattice to a diagonal one, and to carry over all the results for diagonal sublattices to the general case. Indeed, by our isomorphism trick we are able   to extend the definition of the degree of autonomy of an $n$-D system to include nondiagonal sublattices as well. We reduce the problem of the paper to a study of symmetries in Section 5. Section 6 is the core of the paper, and we prove the main result on the existence and uniqueness of the coarsest lattice of definition of an $n$-D system. 

Finally, in Section 7 we show the usefulness of the results of  preceding sections by applying them to a study of the classical system theoretic properties of controllability and autonomy. This application serves to illustrate the main point of the paper, namely that the reduction of an $n$-D system from $\Z^n$ to a sublattice is a reduction in its complexity.

We illustrate our results with a series of examples. Some of the phenomena we describe occur in the case of 1-D systems, and we choose 1-D examples to explain them. Other phenomena occur only in several dimensions, and we choose 2-D examples here, for our intention is to explain them in the simplest, and the most transparent, case possible.

\section{diagonal sublattices of $\Z^n$}
We begin our development by first considering the lattice $\Z^n = \{(x_1,\ldots x_n) ~|~ x_i \in \Z, \forall i \}$ of all points in $\R^n$ with integral coordinates. We denote it by $\mathbb{L}$. For $i = 1, \ldots, n$, let $\sigma_i:\mathbb{L} \rightarrow \mathbb{L}$ map $x = (x_1, \ldots, x_i, \ldots, x_n)$ to $(x_1,\ldots, x_i+1, \ldots, x_n)$; it is the shift operator in the $i$-th direction. A monomial $\sigma^{x'} = \sigma_1^{x'_1}\cdots \sigma_n^{x'_n}$, $x'_i \in \Z$ for all $i$, maps points of $\mathbb{L}$ by composition; thus $\sigma^{x'}(x) = (x_1 + x'_1, \cdots, x_n + x'_n)$.

Let $F$ be any field.
Let $A = F\diff$ be the Laurent polynomial ring generated by these shifts, and their inverses. $A$ is an $F$-algebra, it is the ring of partial difference operators on $\mathbb{L}$.
If $a \in A$ equals $\sum_{j=1}^r c_j \sigma^{z_j}$, where $c_j \neq 0$ and $\sigma^{z_j}$ is a monomial in $A$, then its support is defined to be $\supp(a) = \{z_1, \ldots, z_r\} \subset \mathbb{L}$.

Let $\F_{\mathbb{L}}$ denote the set $F^{\Z^n}$ of all $F$-valued functions on $\mathbb{L}$. The operator $\sigma_i$ induces an action on  $\F_{\mathbb{L}}$ by  mapping $f \in \F_\mathbb{L}$ to $\sigma_if$, where $\sigma_if(x) = f(\sigma_i(x))$. A monomial in $A$ acts by composition, and this action extends to an action of $A$ on $\F_\mathbb{L}$, and  gives it the structure of an $A$-module. The attributes of systems we study take values in $\F_\mathbb{L}$, and we call it the space of signals defined on $\mathbb{L}$.

A linear $n$-D system, by definition, is the kernel of a partial difference operator. Thus, if $P(\sig)$ is an $\ell \times k$ matrix whose entries $p_{ij}$ are in $A$, then it defines a map
\[
\begin{array}{cccc}
P(\sig): & \F_\mathbb{L}^k  & \longrightarrow & \F_\mathbb{L}^\ell\\
& f=(f_1, \ldots, f_k) & \mapsto &  (p_1f, \ldots, p_\ell f) ,
 \end{array}
\]
where the $i$-th row $p_i = (p_{i1}, \ldots, p_{ik})$ of $P(\sig)$ acts on $f$ by $p_if = \sum_{j=1}^kp_{ij}f_j$, $i=1, \ldots, \ell$. The $n$-D system we study is the kernel $\ker_{\F_\mathbb{L}}(P(\sig))$ of this operator. It however depends on the $A$-submodule $P$ of $A^k$ generated by the rows of the matrix  $P(\sig)$, and not on the matrix itself. 
Indeed, the above kernel is isomorphic to $\homo_{A}(A^k/P, ~\F_\mathbb{L})$, the isomorphism given by the map
\[
\begin{array}{ccc}
\ker_{\F_\mathbb{L}}(P(\sig))  & \longrightarrow & \homo_A(A^k/P, ~\F_\mathbb{L})  \\ 
 f=(f_1, \ldots, f_k) & \mapsto & \phi_f   ~,
 \end{array}
\]
where $\phi_f([e_i]) = f_i$, $1 \leqslant i \leqslant k$, and where $[e_1], \ldots ,[e_k]$ denote the images of the standard basis $e_1, \ldots ,e_k$ of $A^k$ in $A^k/P$. 
Hence, we denote this kernel  by $\ker_{\F_\mathbb{L}}(P)$; it is the system defined by $P$ in the signal space $\F_\mathbb{L}$. An element $f \in \ker_{\F_\mathbb{L}}(P)$ is a trajectory of the system. If the rows of $P(\sig)$ are considered to be the laws that govern the system, then a trajectory of the system is a signal that satisfies these laws, namely $P(\sig) f = 0$. This interpretation is foundational in Willems \cite{w}. Clearly, if $P$ is contained in an $A$-submodule $P'$ of $A^k$, then $\ker_{\F_\mathbb{L}}(P') \subset \ker_{\F_\mathbb{L}}(P)$. \\

\noindent Remark 2.1. ~We now give another description of the system $\ker_{\F_\mathbb{L}}(P)$. If the monomial  $\sigma_1^{x_1} \cdots \sigma_n^{x_n}$ in $A$ is identified with $(x_1, \ldots ,x_n)$ in $\mathbb{L}$, then the $F$-algebra $A$ can be identified with the $F$-vector space spanned independently by the points of $\mathbb{L}$.  The signal space $\F_\mathbb{L}$ is then isomorphic to the vector space dual of $A$, namely $\homo_F(A, ~F)$. Given elements $\phi \in \homo_F(A, ~F)$ and $a \in A$, define $a\phi \in \homo_F(A, ~F)$ by $a\phi(a') = \phi(aa')$. This gives $\homo_F(A, ~F)$  the structure of an $A$-module, and the above $F$-isomorphism is an isomorphism of $A$-modules.
Similarly, $\homo_F(A^k/P, ~F)$ is an $A$-module, and it follows that
\[{\homo}_A(A^k/P, ~{\homo}_F(A, ~F))  \simeq {\homo}_F(A^k/P \otimes_A A, ~F) \simeq {\homo}_F(A^k/P, ~F)\]  
by tensor-hom adjunction. The system $\ker_{\F_\mathbb{L}}(P)$ is thus isomorphic to  $\homo_F(A^k/P, ~F)$ as $A$-modules.\\

\noindent Remark 2.2. ~As the functor $\homo_F(-, ~F)$ is exact, so is the functor $\homo_A(-,~\F_\mathbb{L})$ on the category of $A$-modules. This implies that $\F_\mathbb{L}$ is an injective $A$-module. Furthermore $\homo_A(M, ~\F_\mathbb{L}) \neq 0$, if $M \neq 0$, and thus $\F_\mathbb{L}$ is an injective cogenerator.

The injective cogenerator property implies that the assignment $P \rightarrow \ker_{\F_\mathbb{L}}(P)$ is an inclusion reversing  bijective correspondence between $A$-submodules of $A^k$ and $n$-D systems in $\F_\mathbb{L}^k$, see for instance \cite{o,pom,z,stek} (references \cite{pom,stek} deals with systems defined by partial differential equations, but the results that we need here for multidimensional systems are idential in the two cases).\\ 


A sublattice is a $\Z$-submodule (i.e subgroup) of $\mathbb{L}$.  If sublattices $\mathbb{S}$ and $\mathbb{S}'$ satisfy $\mathbb{S} \subset \mathbb{S}'$, then we say that $\mathbb{S}$ is coarser than $\mathbb{S}'$. Given $d = (d_1, \ldots ,d_n) \in \Z^n$ with $d_i \geqslant 0$ for all $i$, define $\mathbb{L}_d$  to be the sublattice $\{(x_1, \ldots, x_n) ~| ~x_i = 0, \pm d_i, \pm 2d_i, \ldots\ ~,\forall i \}$. 
We call $\mathbb{L}_d$ a {\it diagonal} sublattice of $\mathbb{L}$. {\it In this section, as well as in the next Section 3, we consider only diagonal sublattices}. This is no loss of generality, because we show in Section 4 that we can reduce the study of systems on general sublattices to the diagonal case via an isomorphism of $\Z^n$. 

Henceforth we denote the set of nonnegative integers by $\N_+$. We denote the set of those $d$ in $\N_+^n$ with exactly $m$ nonzero entries by $\N^n_m$ ($\N^n_n$ is denoted as usual by $\N^n$). 
 When $d$ is in $\N^n_m$, $\mathbb{L}_d$ is a free $\Z$-module of rank $m$. 
When $d$ is in $\N^n$, we say that the sublattice $\mathbb{L}_d$ is of full rank; otherwise  we call it degenerate.

The ring of difference operators on the sublattice $\mathbb{L}_d$ is the $F$-subalgebra $A_d = F[\sigma_1^{d_1}, \sigma_1^{-d_1}, \ldots ,\sigma_n^{d_n}, \sigma_n^{-d_n}]$ of $A$ generated by monomials corresponding to the points of $\mathbb{L}_d$. Such a subalgebra contains multiplicative inverses of all the nonzero monomials in it, and we call it a {\em Laurent subalgebra} of $A$. 
The support of an element in $A_d$ is clearly contained in the sublattice $\mathbb{L}_d$. The set of $F$-valued functions on the sublattice $\mathbb{L}_d$ is  denoted $\F_{\mathbb{L}_d}$. The ring $A_d$ acts on it by shift, and gives it the structure of an $A_d$-module, with respect to which it is an injective cogenerator.


The $F$-algebra $A$ is a free $A_d$-module. It is
finitely generated  if and only if $d$ is in $\N^n$, and then its rank, denoted by $\rho$, is equal to the product $d_1\cdots d_n$.  
The set $B_d =  \{\sigma_1^{x_1} \cdots \sigma_n^{x_n} ~|~ 0 \leqslant x_i \leqslant d_i - 1, \forall i\}$ is a basis for $A$ as an $A_d$-module. The ring $A$ is then an integral extension of $A_d$.



If $d \in \N^n_m$, and if $I = \{i_1, \ldots , i_m\}$ is the set of indices corresponding to the nonzero entries of $d$, then $A_d \simeq F[\sigma_{i_1}^{d_{i_1}}, \sigma_{i_1}^{-d_{i_1}}, \ldots, \sigma_{i_m}^{d_{i_m}}, \sigma_{i_m}^{-d_{i_m}}]$. 
A free basis for $A$ as an $A_d$ module is now $B_d = \{ \sigma_1^{x_1} \cdots \sigma_n^{x_n} ~|~ 0 \leqslant x_i \leqslant d_i - 1, i \in I; ~x_i \in \Z, i \notin I\}$. 



The inclusion $i:A_d \hookrightarrow A$ of $F$-algebras, is an $A_d$-linear map, and injects $A_d$ into $A$ as a direct summand (corresponding to the element 1 in $B_d$).
Applying $\homo_{F}(-, F)$ to $i$ gives a map $\pi: \F_\mathbb{L} \rightarrow \F_{\mathbb{L}_d}$, which is a surjection. It resticts a function on the lattice $\mathbb{L}$ to the sublattice $\mathbb{L}_d$. This restriction commutes with the action of $A_d$. In other words, $\F_\mathbb{L}$ is an $A_d$-module by restriction of scalars, and $\pi$ is an $A_d$-module map.

Similarly, applying $\homo_{F}(-, ~F)$ to the surjection $A \rightarrow A_d$, gives the inclusion $\F_{\mathbb{L}_d} \hookrightarrow \F_\mathbb{L}$. It maps a  function $f$ on the sublattice $\mathbb{L}_d$, to the function on the lattice $\mathbb{L}$ given by $f$ on the sublattice, and 0 on its non-trivial translates. It is an $A_d$-module map.
  
The aforementioned $A_d$-module structure of $\F_\mathbb{L}$ is simple to describe.

\begin{lemma} Let $d \in \N_+^n$. Then the signal space $\F_\mathbb{L}$ is isomorphic, as an $A_d$-module, to the direct product $\prod_{\sigma^x \in B_d}\F_{\mathbb{L}_d}$ (number of factors equal to the cardinality of $B_d$). Thus, if $d \in \N^n$, i.e. if $\mathbb{L}_d$ is a full rank sublattice, then $\F_L \simeq (\F_{\mathbb{L}_d})^\rho$.
\end{lemma}
\noindent Proof: The lattice $\mathbb{L}$ is the disjoint union of the translates $\{\sigma^x(\mathbb{L}_d) |~\sigma^x \in B_d\}$.  A shift defined by a monomial in the subring $A_d$, leaves each of these translates invariant. Thus, the set of $F$-valued functions on any of these translates is an $A_d$-module, isomorphic to $\F_{\mathbb{L}_d}$.
The isomorphism of the statement is now given by restricting an element in $\F_\mathbb{L}$ to each of these translates of $\mathbb{L}_d$.
In other words, $\F_\mathbb{L} \simeq \homo_F(A, ~F) \simeq \homo_F(\bigoplus_{\sigma^x \in B_d}A_d, ~F) \simeq \prod_{\sigma^x \in B_d} \F_{\mathbb{L}_d}$, as $A_d$-modules.

If $d = (d_1, \ldots, d_n) \in \N^n$, then the cardinality of $B_d$ equals $\rho = d_1\cdots d_n$.
\hspace*{\fill}$\square$\\

\section{contraction and extension}
In this section we introduce the notion of the contraction of an $n$-D  system on $\mathbb{L}$ to a diagonal sublattice, and of the extension of a system on a diagonal sublattice to $\mathbb{L}$. In Section 4, we generalise these notions to the case of general sublattices of $\mathbb{L}$.

Let $d$ be in $\N_+^n$, and $i: A_d^k \hookrightarrow A^k$, be the map defined coordinate wise by $i: A_d \hookrightarrow A$ of Section 2. Given an $A$-submodule $P$ of $A^k$, its contraction to $A_d^k$ is $i^{-1}(P) = P \cap A_d^k$, and is denoted by $P^c$, or by $P^c_d$ for emphasis. It is an $A_d$-submodule of $A_d^k$. 

\begin{definition} Let  $\ker_{\F_\mathbb{L}}(P)$ be the system defined on the lattice $\mathbb{L}$ by the submodule $P \subset A^k$. Its contraction to the sublattice $\mathbb{L}_d$ is the system $\ker_{\F_{\mathbb{L}_d}}(P^c)$ defined by the contraction $P^c$ in $A_d^k$. 
\end{definition}

\begin{proposition} The contracted system $\ker_{\F_{\mathbb{L}_d}}(P^c)$ is the restriction of the elements of $\ker_{\F_\mathbb{L}}(P)$ to the sublattice $\mathbb{L}_d$.
\end{proposition}
\noindent Proof:  
 Consider the inclusion $i: A_d^k/P^c \hookrightarrow A^k/P$ of $A_d$-modules. The functor $\homo_F( -, ~F)$ is exact, hence $\pi: \homo_F(A^k/P, ~F) \rightarrow \homo_F(A_d^k/P^c, ~F)$ is a surjection. 
The proposition now follows from Remark 2.1.
\hspace*{\fill}$\square$\\


\begin{proposition} $($\cite{e}$)$ Let $d \in \N^n$, and let $p \subset A$ be a prime ideal. Then $p^c \subset A_d$ is maximal if and only $p$ is maximal. 

Every prime ideal $q$ of $A_d$ is of the form $p^c$ for some prime $p$ in $A$. There is no containment relation between the primes in $A$ which contract to  $q$. 
\end{proposition}
\noindent Proof: These are consequences of the `going-up' theorem of Cohen-Seidenberg, applicable here as the ring $A$ is integral over $A_d$. \hspace*{\fill}$\square$\\

We note that the Laurent polynomial ring $A$ is the localization of the polynomial ring $F[\sigma_1, \ldots, \sigma_n]$ at the multiplicatively closed set $S = \{(\sigma_1 \cdots \sigma_n)^r | r \geqslant 0\}$ of powers of the product $\sigma_1 \cdots \sigma_n$, and hence that the prime ideals of $A$ are in bijective correspondence with prime ideals of the polynomial ring that do not intersect $S$.\\

\noindent Remark 3.1. (i)~Let $d \in \N^n$ so that $A$ is integral over $A_d$, and let $i \subset A$ be a nonzero ideal. Let $b$ be a nonzero element in it. Let $b^r + a_{r-1}b^{r-1} + \cdots  + a_0 = 0$, $a_i \in A_d$, be an equation of integral dependence for $b$ of smallest degree. As $A$ is an integral domain, it follows that $a_0 \neq 0$. Further, $a_0 \in A_d \cap i = i^c$, and it follows that $i^c$ is not the 0 ideal. Thus 0 is the only ideal of $A$ that contracts to the 0 ideal of $A_d$.

(ii)~ The above is not true if $d \in \N^n_m$. For instance the principal ideal $(1 - \sigma_1)$ in $\C[\sigma_1, \sigma_1^{-1}, \sigma_2, \sigma_2^{-1}]$ contracts to the 0 ideal in $\C[\sigma_2, \sigma_2^{-1}]$. 

(iii)~As for $k \geqslant 2$, consider the submodule generated by $(1, \sigma)$ in $\C[\sigma, \sigma^{-1}]^2$. Its contraction to $\C[\sigma^2, \sigma^{-2}]^2$ equals 0. Thus a nonzero submodule could contract to the 0 submodule even when the sublattice is of full rank.

An immediate corollary is the following result on scalar systems. 
 
\begin{corollary} Let $d \in \N^n$, and let $i$ be an ideal of $A$. Then, 

\noindent (i) the scalar system  $\ker_{\F_\mathbb{L}}(i)$ is a finite dimensional $F$-vector space if and only if the contracted system $\ker_{\F_{\mathbb{L}_d}}(i^c)$ is.

\noindent (ii) $\ker_{\F_\mathbb{L}}(i) = \F_\mathbb{L}$ if and only if $\ker_{\F_{\mathbb{L}_d}}(i^c) = \F_{\mathbb{L}_d}$.
\end{corollary}
\noindent Proof: (i) As $A$ is noetherian, the system $\ker_{\F_\mathbb{L}}(i)$ defined by $i$ is a finite dimensional $F$-vector space if and only if the system $\ker_{\F_\mathbb{L}}(\sqrt{i})$ defined by its radical is, and this is so if and only $\sqrt{i}$ is an intersection of maximal ideals (for instance Corollary 3.9 in \cite{stek}). Statement (i) now follows because $\sqrt{i^c} = (\sqrt{i})^c$. 

\noindent (ii) By Remark 2.2, $\F_\mathbb{L}$, $\F_{\mathbb{L}_d}$ are injective cogenerators as $A$, $A_d$ modules, respectively; hence only the 0 ideal of $A$ or $A_d$, defines the system which is all of $\F_\mathbb{L}$ or $\F_{\mathbb{L}_d}$.  The statement now follows from  Remark 3.1 (i).
\hspace*{\fill}$\square$\\

We return to a systematic study of system theoretic properties of a general $\ker_{\F_\mathbb{L}}(P)$, and its contraction, in Section 7, where we need the following facts.


\begin{lemma} If $P$ is a $p$-primary submodule of $A^k$, then $P^c$ is a $p^c$-primary submodule of $A_d^k$.
\end{lemma}
\noindent Proof: Let $q \in A_d^k \setminus P$ (as $P \subsetneq A^k$, it does not contain $A_d^k$). Suppose $a \in A_d$ is such that $aq \in P^c$. As also $q \in A^k \setminus P$, and as $P$ is a primary submodule of $A^k$, there is an integer $r \geqslant 1$, such that $a^rA^k \subset P$. Then, $a^rA^k \cap A_d^k \subset P \cap A_d^k = P^c$, and thus $P^c$ is primary in $A_d^k$.

The above line also implies that $\sqrt{\ann(A_d^k/P^c)} \subset \sqrt{\ann(A^k/P)} \cap A_d = p^c$. The other inclusion being trivial, it follows that $P^c$ is $p^c$-primary in $A_d^k$.
\hspace*{\fill}$\square$\\

Recollect that a prime ideal $p$ of a commutative ring $R$ is an associated prime of an $R$-module $M$ if it is equal to the annihilator $\ann(x)$ of some $x \in M$; the set of its associated primes is denoted $\ass(M)$. 

\begin{corollary} $P$, a submodule of $A^k$. If $\ass(A^k/P) = \{p_1, \ldots, p_r\}$, then $\ass(A_d^k/P^c) \subset \{p_1^c, \ldots, p_r^c\}$.
\end{corollary}
\noindent Proof: If $P = P_1 \cap \cdots \cap P_r$ is an irredundant primary decompostion in $A^k$, where $P_i$ is the $p_i$-th primary component of $P$, then $P_1^c \cap \cdots \cap P_r^c$ is a not necessarily irredundant primary decomposition of $P^c$. 
\hspace*{\fill}$\square$\\ 

\begin{corollary} (i) If $A^k/P$ is torsion free, then so is $A_d^k/P^c$. \\
(ii) Let $d \in \N^n$. Then, if $A^k/P$ is torsion, so is $A_d^k/P^c$.
\end{corollary} 
\noindent Proof:  (i) To say that $A^k/P$  is torsion free is to say that $P$ is 0-primary, and this implies by the above lemma, that $P^c$ is also 0-primary.

(ii) To say that $A^k/P$ is torsion is to say that 0 is not in $\ass(A^k/P)$. As $d\in \N^n$, the statement now follows from Remark 3.1 (i) and the above corollary. \hspace*{\fill}$\square$\\

We now introduce the notion of extension. If $Q$ is a submodule of $A_d^k$, let $Q^e$ be the $A$-submodule of $A^k$ generated by $Q$. It is the extension of $Q$ to $A^k$, and we describe it explicitly. As $A$ is a free $A_d$-module, tensoring the exact sequence $0 \rightarrow Q \rightarrow A_d^k \rightarrow A_d^k/Q \rightarrow 0$ of $A_d$-modules with $A$, yields the exact sequence 
$ 0 \rightarrow A \otimes_{A_d} Q \rightarrow  A^k \rightarrow A \otimes_{A_d} A_d^k/Q \rightarrow 0$.
Hence $Q^e \simeq A\otimes_{A_d} Q$, and $A^k/Q^e \simeq A \otimes_{A_d} A_d^k/Q$.

\begin{definition} Let  $\ker_{\F_{\mathbb{L}_d}}(Q)$ be the system defined on the sublattice $\mathbb{L}_d$ by the submodule $Q \subset A_d^k$. Its extension to the lattice $\mathbb{L}$ is the system $\ker_{\F_\mathbb{L}}(Q^e)$ defined by the extension $Q^e$ in $A^k$.
 \end{definition}

The extended system admits an elementary description.
 
\begin{proposition} Let $d \in \N_+^n$. The extension $\ker_{\F_\mathbb{L}}(Q^e)$ to the lattice $\mathbb{L}$ is isomorphic to the product $\prod_{\sigma^x \in B_d} \ker_{\F_{\mathbb{L}_d}}(Q)$. 
\end{proposition}
\noindent Proof: The extended system is $\homo_A(A^k/Q^e, ~\F_\mathbb{L})$, and in light of the above observation, is isomorphic to
\[{\homo}_A(A \otimes_{A_d} A_d^k/Q, ~\F_\mathbb{L}) \simeq {\homo}_{A_d}(A_d^k/Q, ~{\homo}_A(A, ~\F_\mathbb{L})) \simeq {\homo}_{A_d}(A_d^k/Q, ~\F_\mathbb{L})\]
by tensor-hom adjunction. By Lemma 2.1, the $A_d$-module structure of $\F_\mathbb{L}$ is isomorphic to the product $\prod_{\sigma^x \in B_d} \F_{\mathbb{L}_d}$, one copy for each of the translates of the sublattice $\mathbb{L}_d$ in $\mathbb{L}$. Hence, $\homo_A(A^k/Q^e, ~\F_\mathbb{L})$ is isomorphic  to  $\prod_{\sigma^x \in B_d} \homo_{A_d}(A_d^k/Q, ~\F_{\mathbb{L}_d}) \simeq \prod_{\sigma^x \in B_d} \ker_{\F_{\mathbb{L}_d}}(Q)$.
\hspace*{\fill}$\square$\\

We illustrate the above proposition with an example that we return to in subsequent sections. \\

\noindent {\bf Example 3.1.} Let $n = 1$ and $d \geqslant 2$; then  $\mathbb{L}_d$ is the sublattice $\{0, \pm d, \pm 2d, \ldots \}$ of $\mathbb{L} = \Z$. The ring of difference operators on  $\mathbb{L}$ (with complex coefficients) is $A = \C[\sigma, \sigma^{-1}]$, and on $\mathbb{L}_d$ is $A_d = \C[\sigma^d, \sigma^{-d}]$. The $F$-algebra $A$ is a free $A_d$-module of rank $d$, and $B_d = \{1, \sigma, \ldots, \sigma^{d-1}\}$ is a basis. The lattice $\mathbb{L}$ is the union of the sublattice $\mathbb{L}_d$ and its translates $\sigma(\mathbb{L}_d), \ldots, \sigma^{d-1}(\mathbb{L}_d)$. 

Let $i \subset A_d$ be the ideal generated by $a = \sigma^d - 1$. A function $f \in \F_{\mathbb{L}_d}$ is a trajectory of the system $\mathcal{I} = \ker_{\F_{\mathbb{L}_d}}(i)$ on $\mathbb{L}_d$ if and only if $(\sigma^d - 1)f = 0$. Thus $\mathcal{I} \simeq \C$, the space of constant functions on $\mathbb{L}_d$.

The extension $i^e$ of $i$ is the ideal generated by $\sigma^d - 1$ in $A$. It defines $\ker_{\F_\mathbb{L}}(i^e)$, the extension of $\ker_{\F_{\mathbb{L}_d}}(i)$ to $\mathbb{L}$, which by Proposition 3.3 is isomorphic to $\prod_{\sigma^j \in B_d} \mathcal{I} \simeq \C^d$. Thus, a trajectory of the extended system is a `piecewise constant' function on $\mathbb{L}$, namely an element of $\F_\mathbb{L}$ which is constant on each translate $\sigma^z(\mathbb{L}_d)$, $z = 0, \ldots, d-1$.

We verify this by a direct calculation. Let $\sigma^d - 1 = (\sigma - \zeta_0) \cdots (\sigma - \zeta_{d-1})$ be the factorization of $a$ in $A$, where $\zeta_0 (=1), \ldots, \zeta_{d-1}$ are the $d$-th roots of unity in $\C$. As $\F_\mathbb{L}$ is an injective $A$-module, it follows that $\ker_{\F_\mathbb{L}}(i^e) = \sum_{i=0}^{d-1} \ker_{\F_{\mathbb{L}}}(\sigma - \zeta_i)$ (see for instance \cite{stek}). A trajectory $f_i$ in $\ker_{\F_\mathbb{L}}(\sigma - \zeta_i)$ satisfies $f_i(x) = \zeta_i f_i(x-1)$, hence $f_i(x) = \zeta_i^x f_i(0)$. As $\zeta_i^d = 1$ for each $i$, it follows that each $f_i$ is piecewise constant, namely the constant $f_i(0)$ on the sublattice $\mathbb{L}_d$, and  the constant $\zeta_i^zf_i(0)$ on the translate $\sigma^z(\mathbb{L}_d)$, $z= 1, \ldots d-1$.
Thus a trajectory of $\ker_{\F_\mathbb{L}}(i^e)$, which is of the general form $f(x) = f_0(0) + \zeta_1^x f_1(0) + \cdots + \zeta_{d-1}^x f_{d-1}(0)$, is also piecewise constant. 

We now show that every piecewise constant function on $\mathbb{L}$ is a trajectory of $\ker_{\F_\mathbb{L}}(i^e)$. Let $c$ be the piecewise constant function which equals $c_z$ on the translate $\sigma^z(\mathbb{L}_d)$, $z = 0, \ldots, d-1$. We require a trajectory $f$ which takes these values $c_z$. In other words, we require $f(z) = f_0(0) + \zeta_1^zf_1(0) + \cdots + \zeta_{d-1}^zf_{d-1}(0) = c_z$, $z = 0, \ldots, d-1$, for appropriate choices of $f_i(0)$. We write this as {\small
\[
\left ( \begin{array}{lccc}
1 & 1 & \cdots & 1 \\ 1 & \zeta_1 & \cdots & \zeta_{d-1}\\ \cdot & \cdot & \cdots & \cdot \\ \cdot & \cdot & \cdots & \cdot \\ 1 & \zeta_1^{d-1} & \cdots & \zeta_{d-1}^{d-1} \end{array} \right ) \left(  \begin{array}{c} f_0(0) \\ f_1(0) \\ \cdot \\ \cdot \\ f_{d-1}(0) \end{array}                     \right) = \left ( \begin{array}{c} c_0 \\ c_1 \\ \cdot \\ \cdot \\ c_{d-1} \end{array} \right )
\]
}
The $d \times d$ matrix  above is invertible, it is the Vandermonde matrix whose determinant equals $\prod_{0 \leqslant i < j \leqslant d-1} (\zeta_i - \zeta_j)$. Hence we can indeed determine the (unique) values of $f_0(0), \cdots, f_{d-1}(0)$ that result in the piecewise constant trajectory $c$. \\

Corresponding to Corollary 3.3 is the next lemma (which follows from $A^k/Q^e \simeq A \otimes_{A_d} A_d^k/Q$). 
  
\begin{lemma} 


\noindent (i) $A_d^k/Q$ is torsion free, or free of rank $s$,  if and only if $A^k/Q^e$ is torsion free, or free of rank $s$ respectively.

\noindent (ii) $A_d^k/Q$ is torsion if and only if $A^k/Q^e$ is torsion.
\end{lemma}


 

Next we study the processes of contraction and extension of a system in conjunction with one another. \begin{lemma} Let $d \in \N_+^n$. Then,

\noindent  (i) for a submodule $P \subset A^k$, $P^{ce} \subset P$, and $P^{cec} = P^c$,

\noindent (ii) for a submodule $Q \subset A_d^k$, $Q = Q^{ec}$.
\end{lemma}
\noindent Proof: (i) is straightforward. As for (ii), it is also straightforward that $Q \subset Q^{ec}$, and that $Q^e = Q^{ece}$. 
 
From the exact sequence $0 \rightarrow Q \rightarrow Q^{ec} \rightarrow Q^{ec}/Q \rightarrow 0$, it follows that $0 \rightarrow A \otimes_{A_d} Q \rightarrow A \otimes _{A_d} Q^{ec} \rightarrow A \otimes _{A_d} Q^{ec}/Q \rightarrow 0$ is also exact. This is just the exact sequence $0 \rightarrow Q^e \rightarrow Q^{ece} \rightarrow A \otimes_{A_d} Q^{ec}/Q \rightarrow 0$. As $Q^e = Q^{ece}$, it follows that $A \otimes _{A_d} Q^{ec}/Q = 0$, and as $A$ is free, it finally follows that $Q = Q^{ec}$. \hspace*{\fill}$\square$\\

The following corollaries are now immediate.

\begin{corollary} Suppose that the submodule $P \subset A^k$ is equal to the extension $Q^e$ of a  submodule $Q \subset A_d^k$. Then $Q = P^c$.
\end{corollary}

\begin{corollary} In Lemma 3.3, $P^{ce} = P$ if and only if $P$ is generated as an $A$-module by elements in $A_d^k$.
\end{corollary}

Let $\mathcal{Q}$ be the system  $\ker_{\F_{\mathbb{L}_d}}(Q)$ on the sublattice $\mathbb{L}_d$, defined by the submodule $Q$ of $A_d^k$. Let $\mathfrak{Q}_\mathbb{L}$ denote the collection of systems on the lattice $\mathbb{L}$ which contract to $\mathcal{Q}$ on $\mathbb{L}_d$. 
As $\F_\mathbb{L}$ is an injective cogenerator, $\mathfrak{Q}_\mathbb{L}$ is in bijective, inclusion reversing, correspondence with the collection $\mathcal{C}(Q) = \{ P \subset  A^k ~ |~ P^c = Q \}$ of $A$-submodules of $A^k$ that contract to $Q$. By Lemma 3,3 (ii), $Q^{ec} = Q$, hence $Q^e \in \mathcal{C}(Q)$. Moreover, if $P \in \mathcal{C}(Q)$, then $Q^e = P^{ce} \subset P$, by Lemma 3.3 (i), and thus $Q^e$ is the unique minimal element of $\mathcal{C}(Q)$. It follows that the system $\ker_{\F_\mathbb{L}}(Q^e)$ is the unique maximal element in $\mathfrak{Q}_\mathbb{L}$. 

As the ring $A$ is Noetherian, the collection $\mathcal{C}(Q)$ contains maximal elements, and correspondingly $\mathfrak{Q}_\mathbb{L}$ contains minimal elements. \\

\noindent {\bf Example 3.2.} Consider the system of Example 3.1.
The ideal $i \subset A_d$ is the principal ideal generated by $a = \sigma^d - 1$. The polynomial $a$ is irreducible in $A_d$, hence $i$ is a maximal ideal. The system $\mathcal{I} = \ker_{\F_{\mathbb{L}_d}}(i)$ on $\mathbb{L}_d$ is isomorphic to $\C$, the space of constant functions on $\mathbb{L}_d$.

As $\sigma^d - 1 = (\sigma - \zeta_0) \cdots (\sigma - \zeta_{d-1})$ is the factorization of $a$ in $A$, where $\zeta_0(=1), \ldots, \zeta_{d-1}$ are the $d$-th roots of unity in $\C$, the collection $\mathcal{C}(i)$ of ideals of $A$ that contract to $i$, are the ideals generated by various products of the $d$ factors of $a$. Thus, the cardinality of $\mathcal{C}(i)$ is $2^d - 1$. The systems defined by these ideals is the collection of systems on $\mathbb{L}$ which contract to $\mathcal{I}$ on $\mathbb{L}_d$.

The extension $i^e$ of $i$ to $A$ defines $\ker_{\F_\mathbb{L}}(i^e)$, the extension of the system $\mathcal{I}$ on $\mathbb{L}$, which by Example 3.1, is isomorphic to $\prod_{\sigma^j \in B_d} \mathcal{I} \simeq \C^d$, the space of `piecewise constant' function on $\mathbb{L}$ (namely an element of $\F_\mathbb{L}$ which is constant on each translate $\sigma^z(\mathbb{L}_d)$, $z = 0, \ldots, d-1$). This system is the unique maximal system on $\mathbb{L}$ which contracts to the system $\mathcal{I}$ on $\mathbb{L}_d$, as $i^e$ is the unique minimal element of $\mathcal{C}(i)$.

We next describe the minimal systems on $\mathbb{L}$ that contract to $\mathcal{I}$ on $\mathbb{L}_d$; these correspond to maximal elements of $\mathcal{C}(i)$, namely the $d$ maximal ideals of $A$ generated by the factors of $a$. Consider the factor $\sigma - \zeta_i, 0 \leqslant i \leqslant d-1$; it defines the system $\ker_{\F_\mathbb{L}}(\sigma - \zeta_i)$ on $\mathbb{L}$. By Example 3.1, a trajectory $f_i$ in it satisfies $f_i(x) = \zeta_i^x f_i(0)$, for all $x \in \mathbb{L}$. We may imagine this trajectory as follows: it is a constant on $\mathbb{L}_d$, equal to $f_i(0)$, and on the translate $\sigma^z(\mathbb{L}_d)$ it equals the constant $\zeta_i^z f_i(0)$. As the point moves from 0 to $d$, the value $f_i(0)$ of $f_i$ at 0, `rotates' till it is back to $f_i(0)$ at $d$. \\

\noindent {\bf Example 3.3.}  Let $d \in \N^n$. By part (i) of Remark 3.1, $0 \subset A$ is the only ideal that contracts to the 0 ideal of $A_d$, hence $\mathcal{C}(0) = \{0\}$. Equivalently, $\F_\mathbb{L}$ is the only system on the lattice $\mathbb{L}$ that contracts to the system $\F_{\mathbb{L}_d}$ on $\mathbb{L}_d$, viz. Corollary 3.1 (ii). Remark 3.1 (ii) and (iii) provide  counter examples when the sublattice is degenerate, or when $k \geqslant 2$.

As for extensions, let $d \in \N_+$.
Then $\F_{\mathbb{L}_d}^k$ is the only system on $\mathbb{L}_d$ that extends to $\F_\mathbb{L}^k$ on $\mathbb{L}$ (as 0 is the only submodule of $A_d^k$ that extends to 0 in $A^k$). \\

We continue our study of $\mathcal{C}(Q)$ in Sections 5 and 7 below.
Now we state the principal result of this section.
\begin{theorem} Suppose that  the system $\ker_{\F_\mathbb{L}}(P)$, defined by the submodule $P \subset A^k$ on the lattice $\mathbb{L}$, is an extension of a system on $\mathbb{L}_d$. Then it can be reconstructed from the contracted  system $\ker_{\F_{\mathbb{L}_d}}(P^c)$ on $\mathbb{L}_d$.
\end{theorem}
\noindent Proof: If $P \subset A^k$ is an extension of a submodule of $A_d^k$, then it is the extension of the submodule $P^c$, by  Corollary 3.4. Hence $\ker_{\F_\mathbb{L}}(P)$ is isomorphic to the product  $\prod_{\sigma^x \in B_d}\ker_{\F_{\mathbb{L}_d}}(P^c)$, by Proposition 3.3.
Thus, $\ker_{\F_\mathbb{L}}(P)$ is determined by the contracted system $\ker_{F_{\mathbb{L}_d}}(P^c)$ on the coarser lattice $\mathbb{L}_d$.  \hspace*{\fill}$\square$\\




\section{nondiagonal sublattices of $\Z^n$}
In this section we consider the case of a general sublattice $\mathbb{S}$ of $\mathbb{L}$ (i.e. not necessarily diagonal). Thus, $\mathbb{S}$ is a free $\Z$-module of rank less than or equal to $n$. If its rank is $m$, then $\mathbb{S}$ is the image of a $\Z$-linear map $S: \Z^m \rightarrow \Z^n$. Let the ring of difference operators on $\mathbb{S}$ be denoted $A_S$. It is the Laurent subalgebra of $A$ generated by the monomials corresponding to points in $\mathbb{S}$. Denote the $F$-valued functions on $\mathbb{S}$ by $\F_\mathbb{S}$; it is an injective $A_S$-module, and also a cogenerator.

In this section we show that $A$ is always a free  $A_S$ module,  for $\mathbb{S}$ an arbitrary sublattice of $\mathbb{L}$, and is finitely generated when $\mathbb{S}$ is of full rank. Indeed, we show that there is a $\Z$-module automorphism of $\mathbb{L}$ which carries the sublattice $\mathbb{S}$ to a diagonal sublattice of $\mathbb{L}$. This allows us to carry over all the results of the previous sections established for diagonal sublattices to the general case. 

\begin{proposition} There is a $\Z$-module automorphism $\phi: \mathbb{L} \rightarrow \mathbb{L}$ which maps the sublattice $\mathbb{S}$ to a diagonal sublattice of ~ $\mathbb{L}$ of the same rank.
\end{proposition}
\noindent Proof: Let the rank of $\mathbb{S}$ be $m$. Let $S: \Z^m \rightarrow \Z^n$ be a $\Z$-linear map whose image is $\mathbb{S}$. Let $S$ also denote the $n \times m$ matrix whose columns are the images of the standard basis for $\Z^m$. Then there exist unimodular matrices $U$ of size $n$ and $V$ of size $m$, with integer entries, and a diagonal $n \times m$ matrix $D$ such that $S = UDV$, namely the Smith Canonical Form of $S$. Let $d_1, \ldots, d_m$ be the diagonal entries of $D$, which we can assume are all strictly positive, as the rank of $S$ equals the rank of $D$. Let $\mathbb{L}_d$ be the diagonal sublattice of $\mathbb{L}$ defined by $d = (d_1, \ldots, d_m)$.

Consider the automorphism $\phi: \mathbb{L} \rightarrow \mathbb{L}$, defined by $\phi(x) = U^{-1}(x)$. An $x \in \mathbb{L}$ is in $\mathbb{S}$ if and only if $x = Sz$ for some $z \in \Z^m$, hence $\phi(\mathbb{S}) = \{U^{-1}Sz~|~ z \in \Z^m \}$, and the latter set is equal to $U^{-1}SV^{-1}V(\Z^m) = DV(\Z^m) \subset \mathbb{L}_d$.

Conversely, $\mathbb{L}_d = D(\Z^m) = U^{-1}SV^{-1}(\Z^m) \subset \phi(\mathbb{S})$. This proves that $\phi(\mathbb{S}) = \mathbb{L}_d$.
\hspace*{\fill}$\square$\\

\begin{corollary} The automorphism $\phi$ of $\mathbb{L}$ in the above proposition induces an $F$-algebra automorphism $\phi_*$ of $A$ which maps $A_S$ to $A_d$. Thus $A$ is a free $A_S$-module, and is finitely generated exactly when $\mathbb{S}$ is a full rank sublattice.
\end{corollary}
\noindent Proof: Define $\phi_*: A \rightarrow A$ by mapping a monomial $\sigma^x$ to $\sigma^{\phi(x)}$, and extending it  linearly to a map of $F$-vector spaces.. It is also an algebra map because $\phi_*(\sigma^x \sigma^{x'}) = \phi_*(\sigma^{x + x'}) = \sigma^{\phi(x + x')} = \sigma^{\phi(x) + \phi(x')} = \sigma^{\phi(x)} \sigma^{\phi(x')} = \phi_*(\sigma^x)\phi_*(\sigma^{x'})$.

As $\phi$ maps points of $\mathbb{S}$ bijectively to points of $\mathbb{L}_d$, it follows that the automorphism $\phi_*$ of $A$ maps  $A_S$ to $A_d$. As $A$ is a free $A_d$-module, it follows that $A$ is also a free $A_S$-module, and is finitely generated only when $\mathbb{S}$ is of full rank.
\hspace*{\fill}$\square$\\

We now extend the map $\phi_*$ to a map on $A^k$, also denoted $\phi_*:A^k \rightarrow A^k$, by defining it coordinate wise: $\phi_*(a_1, \ldots, a_k) = (\phi_*(a_1), \ldots, \phi_*(a_k))$. This is an $F$-vector space map, and is not an $A$-module map. Nonetheless, as $\phi_*: A \rightarrow A$ is an $F$-algebra map, it follows that $\phi_*(a(a_1, \dots, a_k)) =  \phi_*(a)(\phi_*(a_1), \dots, \phi_*(a_k))$. Thus, it follows that if $P$ is an $A$-submodule of $A^k$, its image $\phi_*(P)$ is also an $A$-submodule of $A^k$, and hence that $A^k/P$ is isomorphic to $A^k/\phi_*(P)$ as $F$-vector spaces. By Remark 2.1, $\ker_{\F_\mathbb{L}}(P) \simeq \homo_F(A^k/P, ~F) \simeq \homo_F(A^k/\phi_*(P), ~F) \simeq \ker_{\F_\mathbb{L}}(\phi_*(P))$ as $F$-vector spaces. Furthermore, it also follows that $A^k/P$ is torsion free, or torsion, if and only if $A^k/\phi_*(P)$ is torsion free, or torsion, respectively; indeed, the torsion elements of $A^k/P$ correspond bijectively to the torsion elements of $A^k/\phi_*(P)$ under the above isomorphism. This allows us to carry over all the results of the previous sections about a diagonal sublattice, to the case of a general sublattice. For instance, corresponding to Corollary 3.4, we now have the following statement.
\begin{lemma} Suppose that the submodule $P \subset A^k$ is equal to the extension $Q^e$ of a submodule $Q \subset A_S$. Then $Q = P^c$, the contraction of $P$ to $A_S$.
\end{lemma}

Similarly, we have the following propositions.
\begin{proposition} Let $P \subset A^k$ be an $A$-submodule, and $P^c$ its contraction to $A_S^k$. The contracted system $\ker_{\F_{\mathbb{S}}}(P^c)$ is the restriction of the elements of $\ker_{\F_\mathbb{L}}(P)$ to the sublattice $\mathbb{S}$.
\end{proposition}

\begin{proposition} Let $Q^e$ denote the extension of the submodule $Q$ of $A_S^k$ to $A^k$. Then the system $\ker_{\F_\mathbb{L}}(Q^e)$ on the lattice $\mathbb{L}$ is a direct product of the system $\ker_{\F_\mathbb{S}}(Q)$ on the sublattice $\mathbb{S}$.  
\end{proposition}

\begin{proposition} The system $\ker_{\F_\mathbb{L}}(P)$ is an extension of a system on the sublattice $\mathbb{S}$ if and only if $\ker_{\F_\mathbb{L}}(\phi_*(P))$ is an extension of a system on the diagonal sublattice $\phi(\mathbb{S}) = \mathbb{L}_d$.
\end{proposition}
\noindent Proof: The submodule $P \subset A^k$ is the extension of its contraction to $A_S^k$ if and only if $\phi_*(P)$ is the extension of its contraction to $A_d^k$.    \hspace*{\fill}$\square$\\

The analogues of Lemma 3.1, Corollary 3.3 and Lemma 3.2, with $A_d$ replaced by $A_S$, also continue to hold, an observation we  use in Section 7.\\

We illustrate the above `change of coordinates' by an elementary example. \\

\noindent {\bf Example 4.1.} Let $\mathbb{S}$ be the full rank sublattice of $\Z^2$ generated by $(2,2)$ and $(1, -3)$. It is the image of the map $S: \Z^2 \rightarrow \Z^2$, where  {\small $S = \left ( \begin{array}{lc} 2 & \phantom{x} 1 \\ 2 & -3  \end{array} \right )$}.

We have  {\small $S = \left (\begin{array}{lc} \phantom{C} 1 & 0 \\ -3 & 1 \end{array} \right ) \left ( \begin{array}{lc} 1 & 0 \\ 0 & 8 \end{array} \right ) \left( \begin{array}{lc} 2 & 1 \\ 1 & 0 \end{array} \right )$}. The $\Z$-module automorphism $\phi: \Z^2 \rightarrow \Z^2$ defined by the matrix {\small $\left ( \begin{array}{lc} 1 & 0 \\ 3 & 1 \end{array} \right )$ maps $\mathbb{S}$} to the diagonal sublattice spanned by $(1, 0)$ and$(0, 8)$. The algebra automorphism $\phi_*: A \rightarrow A$ defined by $\phi$, maps $\sigma_1$ to $\sigma_1 \sigma_2^3$, and $\sigma_2$ to itself. It follows that $\phi_*$ maps $A_S = F[\sigma_1^2\sigma_2^2, \sigma_1^{-2}\sigma_2^{-2}, \sigma_1 \sigma_2^{-3}, \sigma_1^{-1}\sigma_2^3]$ to $A_d = F[\sigma_1, \sigma_1^{-1}, \sigma_2^8, \sigma_2^{-8}]$, which is the ring of operators on the diagonal sublattice $\mathbb{L}_d$, where $d = (1, 8)$.

Finally, we check that $A$ is a free $A_S$-module.  A basis for $A$ as an $A_d$-module is $B_d = \{1,\sigma_2, \sigma_2^2, \ldots, \sigma_2^7 \}$. Then $\phi_*^{-1}(B_d) = B_d$ is also a basis for $A$ as an $A_S$-module. \\

We make a final observation about nondiagonal full rank sublattices.

\noindent Remark 4.1.~Let $\mathbb{S}$ be a nondiagonal full rank sublattice of $\mathbb{L}$. Let $S$ be an $n \times n$ integer matrix whose columns are a basis for it. Let its determinant be $\delta$. Then multiplying  it by its adjoint results in the diagonal matrix $\diag(\delta)$. The diagonal sublattice defined by the columns of $\diag(\delta)$, and hence all coarser diagonal sublattices than it, are thus contained in $\mathbb{S}$.

\section{a galois criterion for extended systems}

The purpose of this section is to provide conditions for a system on $\mathbb{L}$ to be an extension of a system defined on a sublattice in terms of symmetries that it must then necessarily  possess. Towards this, we first calculate the group $\aut_{A_d}(A)$ of all $F$-algebra automorphisms of $A$ that fix the Laurent subalgebra $A_d$.

Let $K = F(\sigma) = F(\sigma_1, \ldots, \sigma_n)$ and $K_d = F(\sigma^d) = F(\sigma_1^{d_1}, \ldots, \sigma_n^{d_n})$, denote the fields of fractions of $A$ and $A_d$ respectively, where $d = (d_1, \ldots ,d_n)$. We first study the case when $d \in \N^n$, that is when $\mathbb{L}_d$ is a diagonal sublattice of $\Z^n$ of full rank, postponing the calculation in the general case to the end of the section. Thus, $K$ is now a finite extension of $K_d$, of degree $\rho = d_1 \cdots d_n$. We also confine ourselves to the case of separable extensions.

We denote the $r$-th roots of unity in a ring $R$ by $\mu_r(R)$. 

\begin{proposition} Let $K$ and $K_d$  be as above, where $d \in \N^n$. Let $F$ be either $\C$, or an algebraically closed field of characteristic $p$, with $p$ relatively prime to $d_1, \ldots, d_n$. Then $K$ is a finite Galois extension of $K_d$, whose Galois group is $\gal(K/K_d) = \mu_{d_1}(F) \times \cdots \times \mu_{d_n}(F)$. The group $\aut_{A_d}(A)$ is isomorphic to $\gal(K/K_d)$.
\end{proposition}
\noindent Proof: The irreducible polynomial of $\sigma_i$ over $K_d$ is $(X^{d_i} - \sigma_i^{d_i})$, for each $i$. Hence, as ${\char}(F)$ is relatively prime to $d_i$  (or equals 0), $K$ is separable over $K_d$. As $F$ is algebraically closed, $K$ is also normal over $K_d$, and is thus a finite Galois extension of $K_d$.

Clearly, an automorphism of $K$ which fixes every element of $K_d$, must map $\sigma_i$ to $\zeta_i \sigma_i$, where $\zeta_i$ is a (not necessarily primitive) $d_i$-th root of unity, for each $i$. Thus the Galois group of the extension $K/K_d$ is isomorphic to the group of the statement. It is also clear that each element of $\gal(K/K_d)$ gives by restriction an automorphism of $A$ which fixes $A_d$ pointwise, and hence an element of $\aut_{A_d}(A)$. Conversely, an element of $\aut_{A_d}(A)$ extends to a unique field automorphism of $K/K_d$. Thus $\aut_{A_d}(A) \simeq \gal(K/K_d)$.
\hspace*{\fill}$\square$\\



Let $G$ denote $\aut_{A_d}(A)$. It acts on $A^k$ by
 $g (a_1, \ldots, a_k) = (g(a_1), \ldots, g(a_k))$. It is an $A_d$-module automorphism of $A^k$ that fixes $A_d^k$ pointwise. Let $P$ be an $A$-submodule of $A^k$, then $g(P)$ is also an $A$-submodule of $A^k$. Thus, $G$ acts on the collection of all $A$-submodules of $A^k$ by mapping the submodule $P$ to $g(P)$.

We now provide a group theoretic criterion to determine when a system on $\mathbb{L}$ is the extension of a system on a sublattice $\mathbb{L}_d$ of full rank. We assume henceforth that the field $F$ satisfies the conditions of Proposition 5.1.

\begin{proposition} Let $d \in \N^n$. Let $P$ be an $A$-submodule of $A^k$, and $P^c$ its contraction to $A_d^k$. Then $P = P^{ce}$ if and only if $P$ is $G$-invariant, i.e. $g(P) = P$ for every $g \in G$.
\end{proposition}
We reduce the proof of this proposition to the following lemma.

\begin{lemma}Let $R$ be a (commutative) ring containing all the $r$-th roots of unity, where $r \in \N$ is a unit in $R$. Let $T = R[x, x^{-1}]$ be the Laurent polynomial ring (in one indeterminate) with coefficients in $R$, and let $T_r = R[x^r, x^{-r}]$. Let $i$ be an ideal of $T$, and $i^c$ its contraction to the subring $T_r$. Then $i$ is equal to the extension $i^{ce}$ of $i^c$ to $T$ if $g(i) = i$, for all $g \in \aut_{T_r}(T)$.   
\end{lemma}
\noindent Proof: As in the proof of Proposition 5.1, it follows that $\aut_{T_r}(T) = \mu_{r}(R)$. 
Indeed, an automorphism of $T$ which fixes $T_r$ is given by mapping $x$ to $\zeta_j x$, where $\zeta_j$ is a $r$-th root of unity, $j = 1, \ldots, r$. Denote this automorphism also by $\zeta_j$. 

The ring $T$ is a free $T_r$-module of rank $r$ (with $\{1, x, \ldots, x^{r-1}\}$ a basis). Let $b \in i$, and let $b = b_0 + b_1x + \cdots + b_{r-1}x^{r-1}$, where the $b_j$ belong to $T_r$ for all $j$. By assumption, the elements $\zeta_1(b), \ldots, \zeta_r(b)$ are all in the ideal $i$, and hence so is the sum $\zeta_1(b) + \cdots + \zeta_r(b)$ in $i$. Writing out all the terms in this sum yields 
\[b_0 r+ b_1 (\zeta_1 + \cdots + \zeta_r)x + \cdots + b_{r-1}(\zeta_1^{r-1} + \cdots + \zeta_r^{r-1})x^{r-1} \]
In the above sum, each of the terms in brackets equals 0, as $\zeta_j$ is a $r$-th root of unity, for each $j$. Hence the above sum equals $b_0 r$, and as $r$ is a unit in $R$, it follows that $b_0$ is in the ideal $i$, and hence that it is in $T_r \cap i = i^c$.

It now follows that $b - b_0 = b_1x + \cdots + b_{r-1}x^{r-1}$ is in the ideal $i$, and hence so is $b_1 + \cdots + b_{r-1}x^{r-2}$ also in $i$, as $x$ is a unit in $T$. Repeating the argument in the above paragraph yields that $b_1$ is in the ideal $i$, and hence that it is also in $i^c$. Continuing thus, it follows that $b_0, b_1, \ldots, b_{r-1}$ are all in $i^c$, and hence that $b$ is in $i^{ce}$. This implies that $i = i^{ce}$, as $b$ was an arbitrary element of $i$.  
\hspace*{\fill}$\square$\\ 

\noindent Proof of Proposition 5.2: By Corollary 3.5, $P = P^{ce}$ if and only if $P$ is generated as an $A$-module by elements in $A_d^k$. As $G$ leaves $A_d^k$ fixed pointwise, it follows that $g(P) = P$, for every $g \in G$.
  
We prove the converse by induction on $k$. Let $k = 1$, and let $i$ be an ideal of $A$ such that $g(i) = i$, for all $g \in G$.  Consider the sequence of inclusion of $F$-algebras
\[A_d \hookrightarrow A_d[\sigma_1, \sigma_1^{-1}] \hookrightarrow \cdots \hookrightarrow A_d[\sigma_1, \sigma_1^{-1}, \ldots, \sigma_{n-2}, \sigma_{n-2}^{-1}][\sigma_{n-1}, \sigma_{n-1}^{-1}] \hookrightarrow\]  
\[A_d[\sigma_1, \sigma_1^{-1}, \ldots, \sigma_{n-1}, \sigma_{n-1}^{-1}][\sigma_n, \sigma_n^{-1}] = A . \]
We change notation and rewrite the sequence as
\[ A_d = A_{d,0} \hookrightarrow A_{d,1} \hookrightarrow  \cdots \hookrightarrow A_{d, n-1} \hookrightarrow A_{d,n} = A\]
As $p$ is relatively prime to $d_1, \ldots, d_n$, every inclusion $A_{d,j-1} \hookrightarrow A_{d,j}$, $j=1, \ldots, n$,  satisfies the conditions of Lemma 5.1, with $\aut_{A_{d,j-1}}(A_{d,j}) \simeq \Z/(d_j)$. Then, as the ideal $i \subset A$ is invariant under $\aut_{A_d}(A)$, it is in particular invariant under the subgroup $\aut_{A_{d, n-1}}(A) \simeq \Z/(d_n)$, and hence equals the extension of its contraction to $A_{d, n-1}$. Repeating this argument for every inclusion in the above chain, it finally follows that $i$ is the extension of its contraction $i^c$ to $A_d$.

Now assume that the proposition holds for all $G$-invariant submodules of $A^{k-1}$. Let $P \subset A^k$ be a $G$-invariant submodule of $A^k$. Consider the short exact sequence
\[ 0 \rightarrow A \stackrel{\iota}{\longrightarrow} A^k \stackrel{\pi}{\longrightarrow} A^{k-1} \rightarrow 0
\] 
where $\iota(a) = (a, 0, \ldots, 0)$, and $\pi(a_1, \ldots, a_k) = (a_2, \ldots, a_k)$. Restricting this sequence to $P \subset A^k$ gives the short exact sequence 
\begin{equation} 0 \rightarrow \iota^{-1}(P) \longrightarrow P \longrightarrow \pi(P) \rightarrow 0  
\end{equation}
of $G$-invariant $A$-modules. Similarly, the short exact sequence
\[ 0 \rightarrow A_d \stackrel{\iota}{\longrightarrow} A_d^k \stackrel{\pi}{\longrightarrow} A_d^{k-1} \rightarrow 0
\]
yields the short exact sequence of $A_d$-modules
\[0 \rightarrow \iota^{-1}(P^c) \longrightarrow P^c \longrightarrow \pi(P^c) \rightarrow 0 .
\]
An elementary calculation shows that $\iota^{-1}(P^c) = \iota^{-1}(P)^c$,  $\pi(P^c) = \pi(P)^c$. Extending this exact sequence of $A_d$-modues to $A$-modules by tensoring with $A$, yields the short exact sequence 
\begin{equation} 0 \rightarrow \iota(P)^{ce} \stackrel{\iota}{\longrightarrow} P^{ce} \stackrel{\pi}{\longrightarrow} \pi(P)^{ce} \rightarrow 0 . 
\end{equation}
Sequences (1) and (2) fit into the commutative diagram 
\[
\begin{array}{lcccccccc}
0 & \rightarrow & \iota^{-1}(P) & \longrightarrow & P & \longrightarrow & \pi(P) & \rightarrow & 0 \\
& & \uparrow & & \uparrow & & \uparrow && \\
0 & \rightarrow & \iota^{-1}(P)^{ce} & \longrightarrow & P^{ce} & \longrightarrow & \pi(P)^{ce} & \rightarrow & 0 
\end{array}\]
of $G$-invariant $A$-modules. By induction, $ \iota^{-1}(P)^{ce} = \iota^{-1}(P)$ and $\pi(P)^{ce} = \pi(P)$, hence it follows (by the `Five Lemma') that $P^{ce} = P$.
This concludes the proof of the proposition.
\hspace*{\fill}$\square$\\ 

The above proposition carries over to a general  sublattice $\mathbb{S}$ of full rank. We denote the group of automorphisms of $A$ which fixes every element of $A_S$ (in the notation of Section 4) by $\aut_{A_S}(A)$.

\begin{proposition} Suppose that the sublattice $\mathbb{S} \subset \mathbb{L}$ is of full rank. Let $P$ be a submodule of $A^k$, and let $P^c$ be its restriction to $A_S^k$. Then $P = P^{ce}$ if and only if $P$ is $\aut_{A_S}(A)$-invariant.
\end{proposition} 
\noindent Proof: By Proposition 4.1, $\mathbb{S}$ is isomorphic to a diagonal sublattice $\mathbb{L}_d$, $d \in \N^n$,  via a $\Z$-module automorphism $\phi$ of $\mathbb{L}$. This automorphism induces an $F$-algebra automorphism $\phi_*:A \rightarrow A$, which maps $A_S$ to $A_d$. We define a map $\Phi: \aut_{A_S}(A) \rightarrow \aut_{A_d}(A)$ by $\Phi(g) = \phi_* \circ g \circ \phi_*^{-1}$. It is clearly a group isomorphism.

As observed in Proposition 4.4, $P$ equals the extension of its contraction to $A_S^k$ if and only if $\phi_*(P)$ is the extension of its contraction to $A_d$, which in turn is equivalent to $\aut_{A_d}(A)$-invariance of $\phi_*(P)$. Thus, if $g \in \aut_{A_S}(A)$, it follows that $\Phi(g)(\phi_*(P)) = \phi_* \circ g(P) =  \phi_*(P)$ if and only if $g(P) = P$.
\hspace*{\fill}$\square$\\

We can now state the main result of this section.
\begin{theorem} The system $\ker_{\F_\mathbb{L}}(P)$ on the lattice $\mathbb{L}$ is the extension of a system on the sublattice $\mathbb{S}$ of full rank, if and only if the submodule $P \subset A^k$ is $\aut_{A_S}(A)$-invariant.
\end{theorem}
\noindent Proof: The proof follows  from Lemma 4.1 and Proposition 5.3.
\hspace*{\fill}$\square$\\ 

Recollect from Section 3 the set $\mathcal{C}(Q)$ of $A$-submodules of $A^k$ that contract to the $A_d$-submodule $Q$ of $A_d^k$; they correspond to systems on the lattice $\mathbb{L}$ that contract to the system $\ker_{\F_{\mathbb{L}_d}}(Q)$ on the sublattice $\mathbb{L}_d$.

\begin{lemma} The set $\mathcal{C}(Q)$ is stable under the action of $\aut_{A_d}(A)$. Furthermore, if the sublattice $\mathbb{L}_d$ is of full rank, then the unique fixed point of the action is $Q^e$.
\end{lemma}
\noindent Proof: If $P$ is in $\mathcal{C}(Q)$, then $g(P)^c = g(P^c) = g(Q) = Q$, hence $g(P) \in \mathcal{C}(Q)$, for every $g \in \aut_{A_d}(A)$.
In the case of a full rank sublattice, Proposition 5.2 implies that $P \in \mathcal{C}(Q)$ is a fixed point of the action if and only if $P = P^{ce} = Q^e$.
\hspace*{\fill}$\square$\\

We continue to assume that $\mathbb{L}_d$ is of full rank, i.e. $d \in \N^n$, and now study the action of $\aut_{A_d}(A)$ on $\mathcal{C}(Q) \setminus \{Q^e\}$. Towards this, we establish a `Galois correspondence' result  in which  we use the following terminology:
Let $G = \aut_{A_d}(A) = \mu_{d_1}(F) \times \cdots \times \mu_{d_n}(F)$, as above. Let $H$ be a subgroup of $G$; if $H = H_1 \times \cdots \times H_n$, where $H_i$ is a subgroup of $\mu_{d_i}(F)$ for all $i$, then we say that $H$ is a {\it product} subgroup of $G$

\begin{proposition} Let $d \in \N^n$, and let $G$ be the group $\aut_{A_d}(A)$. There is an inclusion reversing correspondence between subgroups of $G$ and sublattices of $\mathbb{L}$ containing $\mathbb{L}_d$, where product subgroups correspond to diagonal sublattices.
\end{proposition}

\noindent Proof: For a subgroup $H$ of $G$, let $A(H)$ be the subset of $A$ that is fixed by all the elements of $H$; it is an $F$-subalgebra of $A$ that contains $A_d$. As an element of $G$ acts on a monomial term in $A$ by multiplying it by an appropriate element of $F$, it follows that  a Laurent polynomial in $A$ is fixed by an element of $G$ if and only if each monomial term of the polynomial is itself fixed. Thus, $A(H)$ is generated as an $F$-algebra by monomials.  If $\zeta = (\zeta_1, \ldots, \zeta_n)$ is an element of $H$, and $\sigma^x = \sigma_1^{x_1} \cdots \sigma_n^{x_n}$ is a monomial in $A(H)$, then  $\zeta(\sigma^x)\zeta(\sigma^{-x}) = \zeta(\sigma^x \sigma^{-x}) = \zeta(1) = 1$; hence as $\zeta(\sigma^x) = \sigma^x$, it follows that $\zeta(\sigma^{-x}) = \sigma^{-x}$. It follows that $A(H)$ is a Laurent subalgebra of $A$ that contains $A_d$. (Here we have used $\zeta(\sigma^x)$ to denote $(\zeta_1\sigma_1)^{x_1} \cdots (\zeta_n\sigma_n)^{x_n}$.) Furthermore, if $H$ is a product subgroup of $G$, then $A(H) = A_{d'}$ for some $d' = (d_1', \ldots, d_n') \in \N^n$, with $d_i' | d_i$, for all $i$; this follows as $H$ is now generated by elements $(\zeta_1, 1, \ldots 1), (1, \zeta_2, \ldots, 1), \ldots , (1, 1, \ldots, \zeta_n)$, for some $\zeta_i \in \mu_{d_i}$, where the order of $\zeta_i$ equals $d_i'$. Conversely, given a Laurent subalgebra $B$ of $A$ containing $A_d$ that is generated by monomials, the set of elements of $G$ that fix every element of $B$ is a subgroup, say $H$, of $G$. It is clear that $B = A(H)$, and that this is an inclusion reversing bijective correspondence between subroups of $G$, and Laurent subalgebras of $A$ containing $A_d$ that are generated by monomials. In this correspondence, $G$ corresponds to $A_d$ and the identity subgroup corresponds to $A$.

In turn, a Laurent subalgebra $B$ of $A$ containing $A_d$ that is generated by monomials, corresponds to the sublattice $\mathbb{S}$, of full rank, generated by the points of $\mathbb{L}$ corresponding to these monomials, namely the identification of $\sigma_1^{x_1} \cdots \sigma_n^{x_n}$ with $(x_1, \ldots, x_n)$ in Remark 2.1; clearly $\mathbb{S}$ contains $\mathbb{L}_d$. If $\mathbb{S}$ is the image of the map $S: \Z^n \rightarrow \Z^n$, in the notation of Section 4, then $B$ is the ring $A_S$, and if it were the diagonal sublattice $\mathbb{L}_{d'}$, then $B$ is $A_{d'}$. Conversely, a sublattice $\mathbb{S}$ of $\mathbb{L}$ containing $\mathbb{L}_d$, defines the ring $A_S$ of difference operators on it, which is a Laurent subalgebra of $A$ containing $A_d$ that is generated by monomials. This is a bijective correspondence which is inclusion preserving, namely $A_S \subset A_{S'}$ if and  only the sublattices they define satisfy $\mathbb{S} \subset \mathbb{S}'$.

Composing these two correspondences yields an inclusion reversing correspondence between subroups of $\aut_{A_d}(A)$ and sublattices of $\mathbb{L}$ containing $\mathbb{L}_d$, where product subgroups correspond to a diagonal sublattices.  \hspace*{\fill}$\square$\\

We build upon the following example in Section 6. \\

\noindent {\bf Example 5.1.} Let $n=2, ~d=(2,2)$ (thus $\mathbb{L}=\Z^2$, ~$d_1=2, ~d_2=2$). Then $G = \aut_{A_d}(A) = \mu_2(F) \times \mu_2(F)$. The five subgroups of $G$ are: ~$G$, the identity subgroup $\{(1,1)\}$, $\mu_2(F) \times \{1\}$, $\{1\}\times \mu_2(F)$, and the subgroup $H_2 =\{(1,1),(-1,-1)\}$ (we write $\mu_2(F) \simeq \Z/2\Z$ multiplicatively as $\{-1, 1\}$). By the above proposition, the sublattices of $\mathbb{L}$ that contain the sublattice $\mathbb{L}_{(2,2)}$ correspond to the five subgroups, and are, respectively, $\mathbb{L}_{(2,2)}$, $\mathbb{L}$, $\mathbb{L}_{(1,2)}$, $\mathbb{L}_{(2,1)}$ and the sublattice generated by $(1,1)$ and $(2,0)$. This last sublattice corresponding to $H_2$, denoted $\mathbb{H}_2$, is nondiagonal as $H_2$ is not a product subgroup.\\

We return to the description of the set $\mathcal{C}(Q)$ in  Lemma 5.2.

\begin{proposition} Let $P \in \mathcal{C}(Q)$. Let its stabilizer under the action of $G = \aut_{A_d}(A)$ be the subgroup $H$. Then the orbit of $P$ is isomorphic to the quotient $G/H$. When $d \in \N^n$, the orbit of $P$ is finite. In this case, let $\mathbb{H}$ be the sublattice of $\mathbb{L}$ corresponding to the subgroup $H \subset G$. Then $P = P^{ce}$, where $P^c$ is the contraction of $P$ to $A(H)$.
\end{proposition}  
\noindent Proof: The first part is clear. Let $A(H)$ be the ring of difference operators on the sublattice $\mathbb{H}$; then $\aut_{A(H)}(A) = H$, and the result follows from Proposition 5.3.
\hspace*{\fill}$\square$\\

\noindent Remark 5.1.~More generally, we define the set $\mathcal{C}(Q)$  by replacing the diagonal sublattice $\mathbb{L}_d$ with an arbitrary sublattice $\mathbb{S}$. Thus, let $Q \subset A_S^k$, and let $\mathcal{C}(Q) = \{ P \subset  A^k ~ |~ P^c = Q \}$ be the set of $A$-submodules of $A^k$ that contract to the submodule $Q$. As in Section 3, $Q^e$ is the unique minimal element of $\mathcal{C}(Q)$. The $A_S$-module $\F_\mathbb{S}$ is an injective cogenerator, hence the elements in $\mathcal{C}(Q)$ are in bijective, inclusion reversing, correspondence with the collection of systems on $\mathbb{L}$ which contract to $\ker_{\F_\mathbb{S}}(Q)$ on the sublattice $\mathbb{S}$.

All the results on $\mathcal{C}(Q)$ continue to hold in this generality. Thus, Proposition 5.3 implies the following generalisation of Lemma 5.2: {\it The set $\mathcal{C}(Q)$ is stable under the action of $\aut_{A_S}(A)$. Furthermore, if the sublattice $\mathbb{S}$ is of full rank, then the unique fixed point of the action is $Q^e$.} \\

We return to Example 3.1 to elucidate the above results. 

\noindent {\bf Example 5.2.} In the notation of Example 3.1, $a = \sigma^d - 1$ is irreducible in $A_d$, and $(\sigma - \zeta_0) \cdots (\sigma - \zeta_{d-1})$  is its factorization in $A$, where $\zeta_0(=1), \ldots, \zeta_{d-1}$ are the $d$-th roots of unity in $\C$. The ideal $i \subset A_d$ is the principal ideal $(a)$. The set $\mathcal{C}(i)$ of ideals in $A$ that contract to $i$ in $A_d$, is of cardinality $2^d-1$, and consists of the ideals in $A$ generated by the various possible products of the factors of $a$. 

The group $G = \aut_{A_d}(A)$ of automorphisms of $A$ keeping $A_d$ fixed is $\mu_d(\C)$, the group of the $d$-th roots of unity. $G$ acts on $\mathcal{C}(i)$; for instance the element $\zeta \in G$ acts on the principal ideal $(\sigma - \zeta_i) \in \mathcal{C}(i)$ to give the principal ideal $(\zeta\sigma - \zeta_i)$. The fixed point of the full group $G$ is clearly the ideal $i^e$. Thus, only the ideal $i^e = (\sigma^d - 1)$ in $\mathcal{C}(i)$ is the extension of its contraction to $A_d$.

If $d$ is prime, then every element in $\mathcal{C}(i) \setminus \{i\}$ is left invariant only by the trivial subgroup 1. Otherwise, suppose $d' > 1$ divides $d$. Let $H$ be the subgroup of $G$ consisting of the $d'$-th roots of unity, say $\zeta_0, \zeta_{i_1}, \ldots \zeta_{i_{d'-1}}$. Then $H$ leaves invariant the ideal generated by the factor $a' = (\sigma - \zeta_0)(\sigma - \zeta_{i_1}) \cdots (\sigma - \zeta_{i_{d'-1}}) = \sigma^{d'}-1$ of $a$. Thus $A(H)$, the ring of difference operators left fixed by $H$ equals $A_{d'}$. The sublattice of $\mathbb{L}$ corresponding to $H$ is $\mathbb{L}_{d'}$. By Proposition 6.4, the ideal $(a')\subset A$ satisfies $(a') = (a')^{ce}$, where $(a')^c$ is the contraction of $(a)$ to $\mathbb{L}_{d'}$. \\

\noindent Remark 5.2.  ~The principal results of this section are all for sublattices of full rank. If we wish to make a similar study of symmetries of $n$-D systems on degenerate sublattices, then we would need analogues of results of this section for purely transcendental extensions. We make a few brief comments to illustrate the difficulties involved.

Assume now that the sublattice $\mathbb{L}_d$ is degenerate, and suppose that $d \in \N^n_m$. Then the Laurent subalgebra $A_d$ is isomorphic to $F[\sigma_1^{d_1}, \sigma_1^{- d_1}, \ldots, \sigma_m^{d_m}, \sigma_m^{- d_m}]$, and $A$ is not integral over $A_d$. The extension $K = F(\sigma_1, \ldots, \sigma_n)$ of $K_d = F(\sigma_1^{d_1}, \ldots \sigma_m^{d_m})$ is a Galois extension followed by a transcendental extension, and  it suffices for us now to study purely transcendal extensions.

The field of fractions  $K = F(\sigma_1, \ldots, \sigma_n)$ of $A$ is a purely transcendental extension of $F$. The automorphism group $\aut_F(K)$ is the Cremona group, the group of birational transformations of projective space $\mathbb{P}_F^n$. The group $\pgl_{n+1}(F)$ of linear projective transformations is contained in $\aut_F(K)$, and is equal to it only for $n = 1$. In general, $\aut_F(K)$ is an object of current study in algebraic geometry, see for instance \cite{d}. Similarly, the automorphism group of the polynomial ring $F[\sigma_1, \ldots, \sigma_n]$ is also a subject of current research, where the Jacobian Conjecture still remains open \cite{es}.

In contrast, the automorphism group $\aut_F(A)$ of the Laurent polynomial ring $A$, which sits between $ F(\sigma_1, \ldots, \sigma_n)$ and $ F[\sigma_1, \ldots, \sigma_n]$, can be easily determined, as we show next. 

An $F$-algebra endomorphism of $A$ must map $\sigma_i$ to a unit, and hence to a monomial $r_i\sigma_1^{m_{i1}} \cdots \sigma_n^{m_{in}}$, where the $r_i \in F^*$ and the $m_{ij} \in \Z$, for $i, j = 1, \ldots, n$.
Clearly, if this endomorphism is to be an $F$-algebra automorphism, then the $n \times n$ matrix $(m_{ij})$ must be unimodular.
Thus, we can consider this automorphism to be a `composition' of two $F$-algebra automorphisms, the first  a homothety, $\sigma_i \mapsto r_i\sigma_i$, and the second given by $\sigma_i \mapsto \sigma_1^{m_{i1}} \cdots \sigma_n^{m_{in}}$, for all $i$. They define two group homomorphisms:

(i) $\psi_1: (F^*)^n \rightarrow \aut_F(A)$, defined by $\psi_1(R)(\sigma_i) = r_i\sigma_i$, $i = 1, \ldots, n$, where $R = (r_1, \ldots, r_n)$,

(ii) $\psi_2: \gl \rightarrow \aut_F(A)$, defined by $\psi_2(M)(\sigma_i) = \sigma_1^{m_{i1}} \cdots \sigma_n^{m_{in}}$, $i = 1, \ldots, n$, where the entries of $M$ are the $m_{ij}$. 

The two actions they define do {\it not} commute, namely, $\psi_2(M) \circ \psi_1(R) \neq \psi_1(R) \circ \psi_2(M)$, for all $R \in (F^*)^n, M \in \gl$. Hence they do not lift to an action of the product group $(F^*)^n \times \gl$ on $A$. 

There is however an action of a semi-direct product of $(F^*)^n$ and $\gl$ on $A$, defined as follows.

Let $\phi: \gl \rightarrow  \aut((F^*)^n)  $ be the group homomorphism defined by $\phi(M)(R) = (\prod_{i =1}^n r_i^{m_{1i}}, \ldots, \prod_{i=1}^n r_i^{m_{ni}})$, where $R = (r_1,\ldots, r_n)$ and $M$ is the matrix $(m_{ij})$. Let $(F^*)^n \rtimes_\phi \gl$ be the semi-direct product of $(F^*)^n$ and $\gl$ determined by $\phi$.  A routine calculation shows that 
\[ \psi_1(R)\circ \psi_2(M)(\sigma_i) = \psi_2(M) \circ  \psi_1(\phi(M)(R))(\sigma_i), \]
for all $i$, and for all $R \in (F^*)^n, M \in \gl$. Thus, the actions defined by $\psi_1$ and $\psi_2$ lift to the semidirect product, and gives a homomorphism $\Psi: (F^*)^n \rtimes_\phi \gl \rightarrow \aut_F(A)$.

Conversely, an element $g$ of $\aut_F(A)$ is determined by the images of the $\sigma_i$ under $g$. If $g(\sigma_i)$ is equal to $r_i\sigma_1^{m_{i1}} \cdots \sigma_n^{m_{in}}$, for $i = 1, \ldots, n$, then the $g(\sigma_i)$ determine the element $R = (r_1, \ldots, r_n) \in (F^*)^n$ and the unimodular matrix $M = (m_{ij})$. This defines $\aut_F(A) \rightarrow (F^*)^n \rtimes_\phi \gl$, by mapping $g$ to $(R, M)$, which is inverse to $\Psi$.

We have thus established the following result.

\begin{proposition} $\aut_F(A) \simeq (F^*)^n \rtimes_\phi \gl$.
\end{proposition}


Degenerate nondiagonal sublattices of $\Z^n$ are important for several reasons, for instance \cite{mp} discusses restrictions of a system to degenerate sublattices of rank $n$ minus the degree of autonomy of the system. These sublattices play the role of minimal initial conditions required to solve the  difference equations defining the system.

\section{the coarsest lattice of definition}
We now prove the central result of the paper, viz. that there is a coarsest lattice of definition of an $n$-D system.

Towards this, we study the  process of contraction to coarser and coarser sublattices of $\mathbb{L}$, starting with the case of diagonal sublattices.   
Let $I =\{i_1, \ldots, i_r \}$ be a subset of $ \{ 1, \ldots, n \}$, and $I' = \{j_1, \ldots, j_s\}$ its complement, where $r + s = n$. Let $\N^n_I = \{ (d_1, \ldots, d_n)~|~d_i = 1, \forall i \in I, d_i \in \N, \forall  i \in I' \}$. Denote $F[\sigma_{i_1}, \sigma_{i_1}^{-1}, \ldots, \sigma_{i_r}, \sigma_{i_r}^{-1}]$ by $A_I$.



The following lemma is elementary.
\begin{lemma} An $A$-submodule $P$ of $A^k$ is an extension of a submodule of $A_d^k$, for every $d = (d_1, \ldots,d_n) \in \N^n_I$,  if and only if it is an extension of a submodule of  $A^k_I$.
\end{lemma}
\noindent Proof: Again, by Corollary 3.4 it suffices to show that  $P = (P^c_d)^e, \forall d \in \N^n_I$, if and only if $P = (P^c_I)^e$; here $P^c_I$ denotes the restriction of $P$ to $A_I^k$. If the latter holds, then certainly $P = (P^c_d)^e, \forall d \in \N^n_I$.

To prove the converse, it suffices to observe that there is a sequence $d^{(1)}, d^{(2)}, \ldots, d^{(m)} \ldots$ in $\N^n_I$ such that $A_{d^{(1)}}^k \supsetneq A_{d^{(2)}}^k \ldots \supsetneq A_{d^{(m)}}^k \ldots $, and such that the intersection of this nested sequence is $A_I^k$. Thus if $P$ can be generated by elements in every $A_{d^{(m)}}^k$ of the nested sequence, than it can also be generated by elements in their intersection $A_I^k$.  
 \hspace*{\fill}$\square$\\




  
\begin{corollary} 
An $A$-submodule $P$ of $A^k$ is an extension of a submodule of $A_S^k$, for every sublattice $\mathbb{S} \subset \mathbb{L}$,  if and only if it is generated by elements in $F^k$.
\end{corollary}
\noindent Proof:  By Remark 4.1, it suffices to consider the statement for diagonal sublattices of $\mathbb{L}$. Now let $I = \emptyset$ in the above lemma; then $A_\emptyset = F$.     \hspace*{\fill}$\square$\\

\begin{definition} An $A$-submodule $P$ of $A^k$ is said to be constant if it has a set of generators  in $F^k$. The set of constant submodules of $A^k$ is denoted $\mathcal{M}(F)$.
\end{definition}


\begin{lemma} A constant submodule of $A^k$  is free. 
\end{lemma}
\noindent Proof: Let $P$ be a constant submodule, and let $p_1, \ldots, p_\ell$ be a set of generators for $P$ as in the above definition. Without loss of generality, we can assume the $p_i$ to be $F$-independent.

Suppose $a_1p_1 + \ldots, + a_\ell p_\ell = 0$, for some $a_1, \ldots, a_\ell$ in $A$. Write $a_i = c_i + b_i$,  $i = 1, \ldots, \ell$, where $c_i$ is the constant term of $a_i$. Then $\sum_i(c_ip_i + b_ip_i) = 0$ implies $\sum_ic_ip_i = 0$, and thus that $c_i =  0$, for all $i$. Thus none of the $a_i$ has a constant term.

Repeating this argument for each monomial term $\sigma_1^{j_1} \cdots \sigma_n^{j_n}$, proves that the $a_i$ are all 0, and thus that these generators are $A$-independent as well.
\hspace*{\fill}$\square$\\


Let $\mathbb{T}$ and $\mathbb{S}$ be arbitrary sublattices of $\mathbb{L}$.
By Proposition 4.1, we may assume that $\mathbb{T}$ is a diagonal sublattice $\mathbb{L}_d$, for some $d \in \N_+^n$. The corresponding Laurent subalgebras of $A$ are $A_d$ and $A_S$. Let $A_{S_d}$ denote  $A_d \cap A_{S}$, it is the ring of difference operators on the intersection $\mathbb{L}_d \cap \mathbb{S}$, denoted $\mathbb{S}_d$. It is  the largest Laurent subalgebra of $A$ that is contained in both $A_d$ and $A_{S}$. 
 
We expand our previous notation: if $\mathbb{S} \subset \mathbb{S}'$, so that $A_{S} \subset A_{S'}$, then the contraction of $Q \subset A_{S'}^k$ to $A_{S}^k$ is denoted $Q^c_{S}$, or by $Q^c_{d}$ if $\mathbb{S}$ is the diagonal sublattice $\mathbb{L}_{d}$. Similarly, the extension of $R \subset A_{S}^k$ to $A_{S'}^k$ is denoted by $R^e_{S'}$, or by $R^e_{d'}$ if $\mathbb{S}'$ is $\mathbb{L}_{d'}$. We continue to denote by $R^e$ its extension to $A^k$.

\begin{proposition} Let $d = (d_1, \ldots, d_n) \in \N_m^n, ~0 \leqslant m \leqslant n$, and $\mathbb{S}$ be a sublattice of $\mathbb{L}$. Let $i$ be an ideal of $A$. If $i = (i^c_d)^e$, then it follows that $i^c_{S} = (i^c_{S_d})^e_{S}$. 

Thus, if also $i = (i^c_{S})^e$, then $i  = (i^c_{S_d})^e$. 
\end{proposition}
\noindent Proof: The containment $i^c_{S} \supset (i^c_{S_d})^e_{S}$ is trivial, and it remains to show the other containment.

Let $a \in i$; then $a = \sum_j a_jb_j$ for some $a_j \in A, ~b_j \in i^c_d$, as $i = (i^c_d)^e$ by assumption. Now recollect from the discussion preceding Lemma 2.1 that $B_d = \{ \sigma_1^{x_1} \cdots \sigma_n^{x_n} ~|~ 0 \leqslant x_j \leqslant d_j - 1, i \in J; ~x_i \in \Z, j \notin J\}$, where $J$ is the set of indices corresponding to the nonzero entries of $d \in \N_m^n$,
is a basis for $A$ as an $A_d$-module. Thus each $a_j$ in the above sum can be expressed as $a_j = \sum_{\sigma^x \in B_d}  \sigma^xb_x$, where $b_x \in A_d$ (only finitely many of which are nonzero), and $x = (x_1, \ldots, x_n)$, $\sigma^x = \sigma_1^{x_1} \cdots \sigma_n^{x_n}$. Substituting these values of $a_j$ in the first sum and gathering terms, it follows that $a$ can be written as $\sum_{\sigma^x \in B_d} \sigma^xb'_x,~b'_x \in i^c_d$. As the support of each nonzero $b'_x$ is in $\mathbb{L}_d$, it follows that $\supp(\sigma^xb'_x)$ is contained  entirely in the coset $\sigma^x(\mathbb{L}_d)$ of $\mathbb{L}_d$ in $\mathbb{L}$. Thus the support of the distinct terms of $a$ in the above sum are contained in distinct cosets of $\mathbb{L}_d$ in $\mathbb{L}$, and hence the support of $a$ is the disjoint union of these supports.

Now suppose that $a \in  i^c_{S}$, then the support of $a$ is also contained in the sublattice $\mathbb{S}$.  Let 
$z_1, z_2$ be any  two elements in $\supp(\sigma^xb'_x) \subset \sigma^x(\mathbb{L}_d)$. Then $z_1- z_2$ is in $\mathbb{S} \cap \mathbb{L}_d = \mathbb{S}_d$, and hence it follows that the support of the term $\sigma^xb'_x$ of $a$ is now contained in a single coset, say $\sigma^{x'}(\mathbb{S}_d)$, of $\mathbb{S}_d$ in $\mathbb{S}$, for some $x'\in \mathbb{S}$. We rewrite this term as $(\sigma^{x-x'}b'_x)\sigma^{x'}$, so that $\supp(\sigma^{x-x'}b'_x) \subset \mathbb{S}_d$. Hence, $\sigma^{x-x'}b'_x \in i \cap  A_{S_d}$, and as $\sigma^{x'} \in A_{S}$, it follows that this term is in the extension $(i^c_{S_d})^e_{S}$ of $i^c_{S_d}$ to $A_{S}$, and hence so is $a$ in it. 

By Corollary 3.5, the above conclusion is equivalent to saying that the ideal $i^c_{S}$ is generated by elements of its contraction $i^c_{S_d}$ to $A_{S_d}$. If also $i = (i^c_{S})^e$, that is if $i$ is generated by elements in $i^c_{S}$, then it follows that $i$ is generated by elements of $i^c_{S_d}$. This concludes the proof of the proposition.
\hspace*{\fill}$\square$\\

\begin{corollary} Let $P \subset A^k$ be such that both $P = (P^c_d)^e$ and $P = (P^c_{S})^e$. Then $P = (P^c_{S_d})^e$.
\end{corollary}
\noindent Proof: 
As the composition of two restrictions is a restriction, and the composition of two extensions an extension, the proof of the corollary follows by induction on $k$ and the `Five Lemma', exactly as in the proof of Proposition 5.2 above.
\hspace*{\fill}$\square$\\

We can now prove our main result.

\begin{theorem}   If $P$ is in $\mathcal{M}(F)$, then the system $\ker_{\F_\mathbb{L}}(P)$ can be reconstructed from its contraction to any sublattice of $\mathbb{L}$, in particular from its contraction to the sublattice $0 \subset \mathbb{L}$.

If $P$ is not in $\mathcal{M}(F)$, then there is a unique coarsest sublattice $\mathbb{S}$ of $L$, of rank greater than or equal to 1,  such that $\ker_{\F_\mathbb{L}}(P)$ can be reconstructed from its contraction 
$\ker_{\F_{\mathbb{S}}}(P^c_S)$ to $\mathbb{S}$, but not from its contraction to any coarser sublattice of $\mathbb{L}$.
\end{theorem}
\noindent Proof:  The first statement follows from Lemma 6.1 and Theorem 3.1.

Suppose now that $P \notin \mathcal{M}(F)$. Let $\mathfrak{S}$ be the collection of sublattices $\mathbb{S}$ of $\mathbb{L}$ such that the system defined by $P$ can be reconstructed from its contraction to $\mathbb{S}$. Let the minimum of the ranks of all the sublattices in $\mathfrak{S}$ be $r$; then $r \geqslant 1$. Let $\mathbb{S}_1, \mathbb{S}_2, \ldots $ be the sublattices of this minimum rank, listed in some order. Then we can replace $\mathbb{S}_2$ by $\mathbb{S}_1 \cap \mathbb{S}_2$, as by Proposition 4.4 and Corollary 6.2, it is again in $\mathfrak{S}$, and hence also of rank $r$. Similarly, replacing $\mathbb{S}_3$ by $\mathbb{S}_1 \cap \mathbb{S}_2 \cap \mathbb{S}_3$, and so on,  we may assume that there is a nested sequence $\mathbb{S}_1 \supset \mathbb{S}_2 \supset \cdots \supset \mathbb{S}_i \supset \cdots$of sublattices in $\mathfrak{S}$, of lowest rank $r$. By a final application of Proposition 4.1, we may also assume that $\mathbb{S}_1 = \Z^r$, and hence that the nested sequence above is a sequence of full rank sublattices in $\Z^r$. We claim that there are only finitely many terms in this sequence.

By Remark 4.1, we can choose a diagonal sublattice $\mathbb{D}_i \subset \mathbb{S}_i$, for each $i \geqslant 2$. By Proposition 5.4, there are only finitely many sublattices in $\Z^r$ that contain $\mathbb{D}_i$, for any $i$, hence the nested sequence $\Z^r \supset \mathbb{D}_2 \supset \cdots \supset \mathbb{D}_i \supset \cdots$ is cofinal in the original nested sequence. Thus the system can be reconstructed from its contraction to every $\mathbb{D}_i$.
If the sequence of the $\mathbb{D}_i$ has infinitely many terms, then its intersection (or inverse limit) $\mathbb{D}$ is of rank strictly less than $r$, and the system could be reconstructed from  it by Lemma 6.1. This contradicts the minimality of $r$.

Thus there are only finitely many sublattices $\mathbb{S}_i$ of minimum rank $r$ in $\mathfrak{S}$, and their intersection is the unique coarsest sublattice of $\mathbb{L}$ from which the system can be reconstructed.
\hspace*{\fill}$\square$\\
 
The coarsest sublattice from which a system can be reconstructed may of course be coarser than the coarsest diagonal sublattice from which it can be reconstructed.\\

\noindent {\bf Example 6.1.} Consider the scalar system $\ker_{\F_\mathbb{L}}(i)$ on $\mathbb{L} = \Z^2$ defined by the principal  ideal $i = (1+\sigma_1\sigma_2 + \sigma_2^2) \subset A = F[\sigma_1,\sigma_1^{-1}, \sigma_2,\sigma_2^{-1}]$. Let $\mathbb{H}_2$ be the sublattice of Example 5.1, $\mathbb{H}_2 \supset \mathbb{L}_{(2,2)}$, corresponding to the subgroup $H_2$.  If $A_{H_2}$ is the ring of difference operators on $\mathbb{H}_2$, then the ideal $i$ is invariant under $\aut_{A_{H_2}}(A) = H_2$, and hence $i = i^{ce}$, where $i^c$ is the contraction of $i$ to $A_{H_2}$. Thus the system $\ker_{\F_\mathbb{L}}(i)$ can be reconstructed from the sublattice $\mathbb{H}_2$ by Theorem 5.1. It is easy to see that this nondiagonal sublattice of full rank is the coarsest lattice from which the system can be reconstructed.

\vspace{2mm}

\noindent {\bf Example 6.2.} 
Consider the scalar system $\ker_{\F_\mathbb{L}}(j)$ on $\mathbb{L} = \Z^2$ defined by the ideal $j = (1+\sigma_1\sigma_2) \subset A = F[\sigma_1,\sigma_1^{-1}, \sigma_2,\sigma_2^{-1}]$. Let $\mathbb{H}_2$ be as above. Again, the ideal $j$ is invariant under $\aut_{A_{H_2}}(A) = H_2$, and thus the system $\ker_{\F_\mathbb{L}}(j)$ can be reconstructed from the sublattice $\mathbb{H}_2$. 

Now let $d_r=(2^r,2^r)$ and let $H_{2^r}$ denote the subgroup $\{(\zeta, \zeta^{-1}) ~|~ \zeta \in \mu_{2^r}(F)\}$ of $\aut_{A_{d_r}}(A) = \mu_{2^r}(F) \times \mu_{2^r}(F)$ (where $\mu_{2^r}(F)$ denotes the $2^r$-th roots of unity in $F$). Let $\mathbb{H}_{2^r} \supset \mathbb{L}_{(2^r,2^r)}$ be the sublattice corresponding to $H_{2^r}$ given by Proposition 5.4; this sublattice is generated by $(1,1)$ and $(2^r,0)$. If $A_{H_{2^r}}$ is the ring of difference operators on $\mathbb{H}_{2^r}$, then the ideal $j$ is also invariant under $\aut_{A_{H_{2^r}}}(A) = H_{2^r}$; hence $j = j^{ce}$, where $j^c$ is now the contraction of $j$ to $A_{H_{2^r}}$. Thus the system $\ker_{\F_\mathbb{L}}(j)$ can also be reconstructed from the sublattice $\mathbb{H}_{2^r}$. 

In this way, we construct a decreasing nested sequence of sublattices $\mathbb{H}_{2^r}, r \geqslant 1$, of $\mathbb{L}$, from each of which the system $\ker_{\F_\mathbb{L}}(j)$ can be reconstructed. Then the system can also be reconstructed from the intersection of these sublattices (as in the proof of Lemma 6.1). This intersection is the degenerate nondiagonal sublattice of $\Z^2$ of rank 1 generated by $(1,1)$, namely the sublattice $\mathbb{D} = \{ (x, x) \in \Z^2 ~| ~x = 0, \pm1, \pm2, \ldots \}$.

We finally remark that $\aut_F(A) \simeq (F^*)^2 \rtimes_\phi \gltwo$ by Proposition 5.6, and that the subgroup that leaves $j$ invariant contains the infinite subgroup isomorphic to {\small$\{\left (\begin{array}{lc}  \phantom{x} m & m-1 \\ 1-m & 2-m \end{array} \right ) ~|~ m \in \Z \}$}, as well as the infinite set  {\small$\{\left (\begin{array}{lc}  \phantom{x} m & m+1 \\ 1-m & -m \end{array} \right ) ~|~ m \in \Z \}$} (the latter 
includes the transposition $\sigma_1 \mapsto \sigma_2, \sigma_2 \mapsto \sigma_1$). The proof of Lemma 5.1 thus fails here, and so do all subsequent results which rely upon it. We therefore arrive at the (degenerate) coarsest lattice of definition of $\ker_{\F_\mathbb{L}}(j)$, namely $\mathbb{D}$ above, via a sequence of full rank sublattices, for which our results hold.

\section{control of $n$-D systems}
In this section we relate system theoretic properties of an $n$-D system to those of its contraction or extension. We start with the fundamental notion of controllability due to Willems \cite{w}. The following definition of a controllable $n$-D system due to Wood et al.  is patterned after Willems \cite{w} and Rocha \cite{r}; we refer to Zerz's book \cite{z} for details. 

\begin{definition} \cite{wro} An $n$-D system $\mathcal{B}$ is {\it controllable} if there exists a positive real $\rho$ such that for any subsets $T_1, T_2$ of $\Z^n$ with $d(T_1, T_2) > \rho$, $\mathcal{B}(T_1 \cup T_2) = \mathcal{B}(T_1) \times \mathcal{B}(T_2)$.
\end{definition} 
Here $\mathcal{B}(T)$ denotes the restriction of the trajectories of $\mathcal{B}$ to the subset $T \subset \Z^n$, and $d(T_1, T_2) = \min \{d(x_1, x_2) ~| ~x_1 \in T_1, x_2 \in T_2\}$, where $d(x_1, x_2) = \|x_1 - x_2\|_1$, is the $L_1$ norm on $\Z^n$.

The definition posits a solution to a `patching problem', namely,  that if  $T_1$ and $T_2$ are any two subsets of the lattice $\Z^n$ that are separated by distance at least $\rho$, and given any two trajectories $f_1$ and $f_2$ of the system $\mathcal{B}$ restricted to $T_1$ and $T_2$ respectively, then there is a trajectory $f$ in $\mathcal{B}$ that restricts to $f_1$ on $T_1$ and to $f_2$ on $T_2$.

\begin{theorem} \cite{wro} If the system $\mathcal{B}$ is $\ker_{\F_\mathbb{L}}(P)$, $P$ a submodule of $A^k$, then it is controllable if and only if $A^k/P$ is torsion free.
\end{theorem} 

\noindent Remark 7.1. We first explain that the torsion free condition in the above theorem implies a  result which is, a priori, stronger than controllability, namely the existence of an `image representation'.

Let $P(\sigma, \sigma^{-1})$ be an $\ell \times k$ matrix whose $\ell$ rows generate the submodule $P \subset A^k$. Let $R$ be the set of all relations between its $k$ columns. $R$ is an $A$-submodule of $A^k$; clearly it depends only on the submodule $P$, and not the choice of the matrix $P(\sig)$.  Suppose that $R$ is generated by $k_1$ elements. Let $R(\sigma, \sigma^{-1})$ be the $k \times k_1$ matrix whose columns are these generators. Then the sequence 
\[ A^\ell \stackrel{P(\sig)^\t}{\longrightarrow} A^k \stackrel{R(\sig)^\t}{\longrightarrow} A^{k_1} \]
is a complex (the superscript $\t$ denotes transpose), which is exact if and only if $A^k/P$ is torsion free, \cite{o,stek}.  

Assuming that $A^k/P$ is torsion free, and applying the functor $\homo_A(-,~\F_\mathbb{L})$ to this exact sequence gives the exact sequence
\[ \F_\mathbb{L}^{k_1} \stackrel{R(\sig)}{\longrightarrow} \F_\mathbb{L}^k \stackrel{P(\sig)}{\longrightarrow} \F_\mathbb{L}^\ell ~,  \]
as $\F_\mathbb{L}$ is an injective $A$-module.  Hence $\ker_{\F_\mathbb{L}}(P)$ is equal to the image of $R(\sig)$, and it is elementary that the existence of such an image representation implies controllability. Thus torsion freeness of $A^k/P$ implies that $\ker_{\F_\mathbb{L}}(P)$ admits an image representation, and hence that it is controllable. 

Furthermore, as $\F_\mathbb{L}$ is also a cogenerator, namely Remark 2.2, it follows that if $A^k/P$ is not torsion free, then it is not controllable, and hence that it does not admit an image representation, see \cite{stek}. 
\\

In what follows, we use the notation of Section 4. Thus $\mathbb{S}$ is a sublattice of $\mathbb{L}$, and $A_S$ is the ring of difference operators on it.

\begin{proposition}  If $\ker_{\F_\mathbb{L}}(P)$ is controllable, then the contracted system $\ker_{\F_\mathbb{S}}(P^c)$ on $\mathbb{S}$ is also controllable (where $P^c$ is the contraction of the submodule $P \subset A^k$ to $A_S^k$).  
\end{proposition}
\noindent Proof: We have observed in Section 4 that the analogue of Corollary 3.3 (i) holds for the sublattice $\mathbb{S}$; hence $A_\mathbb{S}^k/P^c$ is torsion free if $A^k/P$ is torsion free. The $A_S$-module $\F_\mathbb{S}$ is an injective cogenerator,  thus the system  $\ker_{\F_\mathbb{S}}(P^c)$ on $\mathbb{S}$ admits an image representation  if $\ker_{\F_\mathbb{L}}(P)$  admits an image representation on $\mathbb{L}$.
\hspace*{\fill}$\square$\\ 

\noindent Remark 7.2. ~The submodule $R \subset A^k$ (in the notation of Remark 7.1) which determines the image representation of  the system $\ker_{\F_\mathbb{L}}(P)$ does not necessarily contract to the submodule which determines the image representation of the contracted system. For instance, consider the example of a nonzero submodule $P \subset A^k$ which contracts to the 0 submodule of $A_S^k$, as in Remark 3.1 (iii). Then $R$ is strictly contained in $A^k$, and hence does not contract to $A_S^k$, which is the submodule that determines the image representation of the contracted system $\ker_{\F_\mathbb{S}}(0)$.

Image representations are however well behaved under extensions, as we discuss below.  

\noindent Remark 7.3.~We briefly discuss systems defined by constant submodules of $A^k$, namely the definition of $\mathcal{M}(F)$ in 6.1. By Lemma 6.2, every $P \in \mathcal{M}(F)$ is free, and in fact satisfies the conditions of Theorem 3.1 of \cite{sh} (the exposition there is for distributed systems defined by constant coefficient partial differential equations,  but the identical results also hold for $n$-D systems). Thus, not only is $\ker_{\F_\mathbb{L}}(P)$ controllable for every $P \in \mathcal{M}(P)$, but it also satisfies a generalisation of the classical Popov-Belevitch-Hautus test \cite{l}.

If $P \in A^k$ is a constant submodule, then by definition, its restriction $P^c$ to the zero sublattice of $\mathbb{L}$ satisfies $P^{ce} = P$, and hence the conclusions of the above paragraph hold on every sublattice of $\mathbb{L}$. 

\begin{definition} \cite{w, wro} An $n$-D system $\mathcal{B}$ is autonomous if no nonzero subsystem of $\mathcal{B}$ is controllable.  

An autonomous system is also called uncontrollable, because such a system does not  admit any inputs \cite{stek}.
\end{definition}
 
\begin{theorem} \cite{wro} If the system $\mathcal{B}$ is $\ker_{\F_\mathbb{L}}(P)$, $P$ a submodule of $A^k$, then it is autonomous if and only if $A^k/P$ is torsion.
\end{theorem} 

\begin{proposition}  Let $\mathbb{S}$ be a sublattice of full rank. If $\ker_{\F_\mathbb{L}}(P)$ is autonomous, then the contracted system $\ker_{\F_{\mathbb{S}}}(P^c)$ on $\mathbb{S}$ is also autonomous.
\end{proposition}
\noindent Proof: The analogue of Corollary 3.3 (ii) holds for the sublattice $\mathbb{S}$ of full rank; hence $A_\mathbb{S}^k/P^c$ is torsion if $A^k/P$ is torsion. \hspace*{\fill}$\square$\\

More generally, for $P$ a submodule of $A^k$, let $P_0$ be the submodule $\{x \in A^k ~| ~\exists ~a \neq 0 ~with ~ax \in P\}$. Then $P_0$ contains $P$, and the quotient $P_0/P$ is the submodule of $A^k/P$ consisting of its torsion elements. The following sequence 
\begin{equation}
0 \rightarrow P_0/P \longrightarrow A^k/P \longrightarrow A^k/P_0 \rightarrow 0
\end{equation}
is exact, where $A^k/P_0$ is torsion free. In general, given a short exact sequence of $A$-modules, the associated primes of the middle term is contained in the union of the associated primes of the other two modules. However, here it is clear that  we have equality.

\begin{lemma} $\ass (A^k/P) = \ass (P_0/P) \sqcup \ass (A^k/P_0)$ (disjoint union). Hence, if $P \subsetneq P_0 \subsetneq A^k$, then $\ass (A^k/P_0) = \{0\}$ and $\ass (P_0/P)$ is the set of all the nonzero associated primes of $A^k/P$.
\end{lemma}

Applying the exact functor $\homo_A(-, ~\F_\mathbb{L})$ to the above sequence gives the exact sequence
\begin{equation}
0 \rightarrow {\ker}_{\F_\mathbb{L}}(P_0) \longrightarrow {\ker}_{\F_\mathbb{L}}(P) \longrightarrow {\homo}_A (P_0/P,~\F_\mathbb{L}) \rightarrow 0 
\end{equation}
By the above lemma $A^k/P_0$ is torsion free, hence $\ker_{\F_\mathbb{L}}(P_0)$ is a controllable sub-system of $\ker_{\F_\mathbb{L}}(P)$. If $P_1$ is any $A$-submodule of $A^k$ such that $P \subset P_1 \nsupseteq P_0$, then $A^k/P_1$  has torsion elements and $\ker_{\F_\mathbb{L}}(P_1)$ is not controllable. Hence $\ker_{\F_\mathbb{L}}(P_0)$ is the largest controllable sub-system of $\ker_{\F_\mathbb{L}}(P)$  in the sense that any other controllable sub-system is contained in it. It is the controllable part of the system $\ker_{\F_\mathbb{L}}(P)$  (see \cite{stek} for more details).

Suppose $P_0/P$ can be generated by $r$ elements; then $P_0/P \simeq A^r/R$, for some submodule $R \subset A^r$, hence $\homo_A(P_0/P, ~\F_\mathbb{L}) \simeq \ker_{\F_\mathbb{L}}(R)$.  As $A^r/R$ is a torsion module, $\ker_{\F_\mathbb{L}}(R)$ is autonomous (this system is a quotient, and not a sub-system, of $\ker_{\F_\mathbb{L}}(P)$, unless the above short exact sequence splits). The sequence (4) is the `controllable-uncontrollable decomposition' of the system $\ker_{\F_\mathbb{L}}(P)$, \cite{wro,stek}. (References \cite{bv,dn,zl} study weaker notions of this decomposition.)\\

We assume now that the sublattice $\mathbb{S}$ is of full rank. Contracting $P \subsetneq P_0 \subsetneq A^k$ to $A_S^k$ gives $P^c \subsetneq P_0^c \subsetneq A_S^k$, where the first strict inclusion is because the quotient $P_0^c/P^c \simeq A_S^r/R^c$ is torsion (by Corollary 3.3), and the second strict inclusion is obvious. Hence, contracting the exact sequence (3) to $A_S^k$ results in the exact sequence
\[ 0 \rightarrow P^c_0/P^c \longrightarrow A_S^k/P^c \longrightarrow A_S^k/P^c_0 \rightarrow 0
\]
where $ A_S^k/P^c_0$ is torsion free, and $P^c_0/P^c$ is the torsion submodule of $A _S^k/P^c$. Thus, applying the exact functor $\homo_{A_d}(-, ~\F_\mathbb{S})$  gives
\[
0 \rightarrow {\ker}_{\F_\mathbb{S}}(P_0^c) \longrightarrow {\ker}_{\F_\mathbb{S}}(P^c) \longrightarrow {\homo}_{A_S} (P_0^c/P^c,~\F_\mathbb{S}) \rightarrow 0
\]
which is the controllable-uncontrollable decomposition of $\ker_{\F_\mathbb{S}}(P^c)$.\\


We now study controllability and autonomy of extensions.

\begin{proposition} Let $Q \subset A_S^k$ be a submodule. The system $\ker_{\F_\mathbb{S}}(Q)$ on the sublattice $\mathbb{S}$ is controllable if and only if its extension  $\ker_{\F_\mathbb{L}}(Q^e)$ to $\mathbb{L}$ is controllable, and is autonomous if and only if the extension is autonomous.
\end{proposition} 
\noindent Proof: We have already observed in Section 4 that Lemma 3.2 holds with $A_d$ replaced by $A_S$; thus $A_S^k/Q$ is torsion free, or torsion, if and only if $A^k/Q^e$ is torsion free, or torsion, respectively.  \hspace*{\fill}$\square$\\

More generally, if 
$0 \rightarrow Q_0/Q \longrightarrow A_S^k/Q \longrightarrow A_S^k/Q_0 \rightarrow 0$
is the analogue of the exact sequence (3), where $Q_0/Q$ is the torsion submodule of $A_S^k/Q$, and $A_S^k/Q_0$ is torsion free, then its extension 
$0 \rightarrow Q_0^e/Q^e \longrightarrow A^k/Q^e \longrightarrow A^k/Q_0^e \rightarrow 0$
is exact, where $Q_0^e/Q^e$ is torsion and $A^k/Q_0^e$ is torsion free. Thus the controllable-uncontrollable decomposition of $\ker_{\F_\mathbb{S}}(Q)$ determines the controllable-uncontrollable decomposition of its extension $\ker_{\F_\mathbb{L}}(Q^e)$.\\

We briefly discuss image representations for extended systems. As explained in Remark 7.1, an image representation for a controllable system defined by a submodule $P \subset A^k$ is determined by a submodule $R \subset A^k$, which is generated by the columns of a matrix $R(\sigma, \sigma^{-1})$ that is a `right annihilator' for $P$. This implies that if $p = (p_1, \ldots, p_k)$ is an arbitrary element of $P$, and $r = (r_1, \ldots, r_k)$ an arbitrary element of $R$, then $p \cdot r = \sum_{i=1}^k p_i r_i = 0$. We write this symbolically as $P \cdot R = 0$.

Now let $Q \subset A_S^k$ define a system on the sublattice $\mathbb{S} \subset \mathbb{L}$. As in Remark 5.1, the submodules of $A^k$ in $\mathcal{C}(Q)$ define $n$-D systems on $\mathbb{L}$ that contract to $\ker_{\F_\mathbb{S}}(Q)$ on $\mathbb{S}$. If $P \in \mathcal{C}(Q)$ defines a controllable system, and if the submodule $R$ defines an image representation for it, then  as $P \cdot R = 0$, it follows that $g(P) \cdot g(R) = 0$, for every $g \in \aut_{A_S}(A)$, and hence that $g(R)$ is the right annihilator for $g(P)$. This implies that $g(R)$ defines an image representation for the controllable system defined by $g(P)$. Thus, controllability and image representations behave well along an  $\aut_{A_S}(A)$-orbit in $\mathcal{C}(Q)$. 

Now assume that $\mathbb{S}$ is of full rank, and suppose that $Q$ defines a controllable system on it. Then $Q^e$ defines a controllable system on $\mathbb{L}$. Let its image representation be defined by the submodule $R$, so that $Q^e \cdot R = 0$. It now follows that $g(Q^e) \cdot g(R) = Q^e \cdot g(R) = 0$ for all $g \in \aut_{A_S}(A)$. This implies that $g(R) = R$ for all $g$, and hence that $R = R^{ce}$, by Proposition 5.3. We have thus established the following proposition.

\begin{proposition} Let $\mathbb{S}$ be a sublattice  of full rank. The submodule $Q \subset A_S^k$ defines a controllable system on $\mathbb{S}$ if and only if its extension to $\mathbb{L}$ defined by $Q^e$ is controllable. Then $T \subset A_S^k$ defines an image representation for $\ker_{\F_\mathbb{S}}(Q)$ if and only if $T^e$ defines an image representation for $\ker_{\F_\mathbb{L}}(Q^e)$.
\end{proposition}

Observe that while controllability is well behaved with respect to arbitrary contractions, viz. Proposition 7.1, autonomy is well behaved only with respect to contractions to full rank sublattices, viz. Proposition 7.2. Indeed, an autonomous $n$-D system on $\mathbb{L}$ may contract to a nonautonomous system on a degenerate sublattice. This is an important phenomenon, and we study it next.

If an $n$-D system on $\mathbb{L}$ admits `inputs', i.e. if it is not autonomous, then it follows from Proposition 7.1 and the exact sequence (4) that its restriction to any degenerate sublattice also admits inputs. On the other hand, a system which is autonomous on $\mathbb{L}$ might become non-autonomous upon restriction to a degenerate sublattice - this phenomenon does not arise for restrictions to sublattices of full rank by Proposition 7.2. This prompts the following definition.

\begin{definition} \cite{wro2,dp} The degree of autonomy of an $n$-D system is the co-rank of the largest diagonal sublattice of $\mathbb{L}$ such that the restriction of the system to it is not autonomous. 
\end{definition}
Thus, the degree of autonomy of a nonautonomous system equals 0, and varies between 1 and $n$ for a nonzero autonomous system (a system is said to be strongly autonomous when the degree equals $n$, \cite{lz, stek}).

We ask the following question: Let $\mathcal{B}$ be an $n$-D system on $\mathbb{L}$, and let $\mathcal{B}^c$ be its restriction to a full rank sublattice $\mathbb{L}_d$. What is the relationship between the degrees of autonomy of $\mathcal{B}$ and $\mathcal{B}^c$~? 

\begin{proposition} The degrees of autonomy of $\mathcal{B}$ and of its restriction $\mathcal{B}^c$ to  the full rank sublattice $\mathbb{L}_d$, are equal.
\end{proposition}
\noindent Proof: If $\mathcal{B}$ is not autonomous, then neither is the restriction $\mathcal{B}^c$ to $\mathbb{L}_d$, as we observed above. Hence, both the degrees of autonomy are equal to 0 in this case. 

Suppose now that $\mathcal{B}$ is autonomous; this implies that $\mathcal{B}^c$ is also autonomous. Suppose that the degree of autonomy of $\mathcal{B}$ equals $r = n-m$, and let $\mathbb{L}_{d'}, ~d' \in \N^n_m$ be a diagonal sublattice of rank $m$ on which the restriction of $\mathcal{B}$ first becomes non-autonomous. Denote this restriction by $\mathcal{B}_m$. The sublattice $\mathbb{L}_{d'} \cap \mathbb{L}_d$ is of full rank in $\mathbb{L}_{d'}$, hence the restriction $(\mathcal{B}_m)^c$ of $\mathcal{B}_m$ to it remains non-autonomous. But $(\mathcal{B}_m)^c$ is also the restriction of $\mathcal{B}^c$ to $ \mathbb{L}_{d'} \cap \mathbb{L}_d$, and hence the degree of autonomy of $\mathcal{B}^c$ is not larger that $r$.

Conversely, if the degree of autonomy of $\mathcal{B}^c$ were smaller than $r$, then there is a sublattice of $\mathbb{L}_d$ of rank $m' > m$ on which its restriction has already become non-autonomous. However, this would imply that $\mathcal{B}$ has also become non-autonomous on a rank $m'$ sublattice of $\mathbb{L}$. This contradiction proves the proposition. 
\hspace*{\fill}$\square$\\

\noindent Remark 7.4. ~The definition of degree of autonomy as well as all the results in \cite{dp,sr,wro2}, and in other papers in the subject, would hold if diagonal sublattices were replaced by arbitrary sublattices. Thus, for instance, the definition of degree of autonomy above could be replaced by the statement `{\em the degree of autonomy of an $n$-D system is the co-rank of the largest sublattice of $\mathbb{L}$ such that the restriction of the system to it is not autonomous}. This is because of the results of Section 4. Using the notation there in the context of the above proposition, if $\phi: \mathbb{L} \rightarrow \mathbb{L}$ is the isomorphism that carries a nondiagonal sublattice $\mathbb{S}$ to the diagonal sublattice $\mathbb{L}_{d'}$, then the induced map $A^k/P \rightarrow A^k/\phi_*(P)$ preserves the properties of controllability and autonomy, and hence also preserves the degree of autonomy of the two systems defined by $P$ and by $\phi_*(P)$. 

\section{Concluding Remarks}
In summary, we have established that an $n$-D system on $\Z^n$ could arise as an extension of a system on a sublattice of $\Z^n$, and that this possibility is equivalent to the existence of a group of symmetries that leave invariant the equations defining the system. These symmetries can be expressed concretely as a subgroup of a Galois group. Such an extension from a sublattice facilitates the study of properties of the $n$-D system, in particular its decomposition into its controllable and autonomous parts. 

We have already pointed out the implications of the results of this paper to the problem of reducing the order of an $n$-D system. These questions, as well as the connections between groups of symmetries and efficient Gr\"{o}bner bases algorithms, will be pursued elsewhere.

An important question that arises here is whether there are analogues of the results of this paper to distributed systems defined by partial differential equations on $\R^n$. There are several problems that would immediately arise in this setting, for instance there is no analogue in $\R^n$ of a full rank sublattice of $\Z^n$. If instead we were to consider a subspace of $\R^n$, then the restriction of the distributed system to this subspace might not be a system \cite{jcw}. Thus, there is no analogue of the notion of the degree of autonomy of an $n$-D system to the distributed case. One way to circumvent these problems might be to directly address these questions at the level of rings of differential operators in the setting of fractional partial differential equations, as alluded to in the introduction, but this would be a research proposal in itself.

\section {Acknowldegement} We are very grateful to Ananth Shankar for his generous help with this paper.
The first author  acknowledges support from the `MATRICS' Grant of the Science and Engineering Research Board, Govt. of India (Project File No. MTR/2019/000907). The second author is grateful to the Department of Electical Engineering for its hospitality during many visits.

\end{document}